\documentclass[12pt,reqno]{article}
 \usepackage[body={6.5in,9.0in},left=1in,top=1in,centering]{geometry}
\usepackage{graphicx,color}
\usepackage{amsmath,amsthm,amssymb,amsfonts,multirow,setspace,epsfig,latexsym}
\usepackage[letterpaper=true,colorlinks=true,urlcolor=blue,linkcolor=blue,citecolor=blue,pdfstartview=FitH]{hyperref}
\usepackage{pstricks}
\usepackage{epstopdf}
\usepackage{mathrsfs}
 \usepackage[sort,comma]{natbib}
\allowdisplaybreaks
\newcmykcolor{orange}{0 .61 0.87 0}
\newcmykcolor{navyblue}{0.94 0.74 0 0}
\newcmykcolor{peach}{0 0.50 0.70 0}

\begin{document}

\makeatletter \@addtoreset{equation}{section}
\renewcommand{\thesection}{\arabic{section}}
\renewcommand{\theequation}{\thesection.\arabic{equation}}

\newcommand{\bdd}{\hspace*{-0.08in}{\bf.}\hspace*{0.05in}}
\def\para#1{\vskip 0.4\baselineskip\noindent{\bf #1}}
\def\qed{\hfill$\qquad \Box$}
\def\rr{{\mathbb R}}\newcommand{\R}{\mathbb R}
\def\ss{{\mathbb S}}  \newcommand{\M}{\mathbb S}

\def\zz{{\Bbb Z}}
\def\G{{\mathcal G}}
\def\fubao#1 {\fbox {\footnote {\ }}\ \footnotetext {From Fubao: #1}}
\def\chao#1 {\fbox {\footnote {\ }}\ \footnotetext {From Chao: #1}}

\newcommand{\bed}{\begin{displaymath}}
\newcommand{\eed}{\end{displaymath}}
\newcommand{\bea}{\bed\begin{array}{rl}}
\newcommand{\eea}{\end{array}\eed}
\newcommand{\disp}{\displaystyle}
\newcommand{\ad}{&\!\!\!\disp}
\newcommand{\aad}{&\disp}
\newcommand{\barray}{\begin{array}{ll}}
\newcommand{\earray}{\end{array}}
\newcommand{\beq}[1]{\begin{equation} \label{#1}}
\newcommand{\eeq}{\end{equation}}
\newcommand{\diag}{\mathrm{diag}}
\newcommand{\La}{\Lambda}
\newcommand{\la}{\lambda}
\newcommand{\sg}{\sigma}
\newcommand{\lf}{\lfloor}
\newcommand{\rf}{\rfloor}
\newcommand{\wdt}{\widetilde}
\newcommand{\wdh}{\widehat}
\newcommand{\bQ}{{\mathbb Q}} \newcommand{\Q}{\mathbb Q}
\newcommand{\al}{\alpha}
\newcommand{\lbar}{\overline}
\newcommand{\tr}{\mathrm {tr}}
\newcommand{\sgla}{\sigma_{\lambda_{0}}}
\newcommand{\lan}{\langle}
\newcommand{\ran}{\rangle}
\newcommand{\comment}[1]{}
\newcommand{\im}{\mathtt{i}}
\newcommand{\E}{{\mathbb E}}
\newcommand{\bP}{{\mathbb P}}   \renewcommand{\P}{\mathbb P}
\newcommand{\m}{\mathfrak m}
\newcommand{\p}{\mathfrak p}
\newcommand{\mz}{{m_{0}}}
\newcommand{\TV}{\mathrm{Var}}

\newcommand{\LL}{{\mathcal L}}
\newcommand{\C}{{\mathcal C}}
\newcommand{\A}{{\mathcal A}}
\newcommand{\B}{{\mathcal B}}
\newcommand{\F}{{\mathcal F}}
\newcommand{\op}{\mathscr A}
\newcommand{\N}{{\mathcal N}} \newcommand{\bN}{{\mathbb N}}
\newcommand{\e}{\varepsilon}
\newcommand{\one}{\mathbf 1}\newcommand{\1}{\mathbf 1}
\renewcommand{\d}{\mathrm{d}}
\newcommand{\set}[1]{\left\{#1\right\}}
\newcommand{\abs}[1]{\left\vert#1\right\vert}
\newcommand{\norm}[1]{\left\|#1\right\|}
\newcommand{\sgn}{\mathrm{sign}}

\newtheorem{Theorem}{Theorem}[section]
\newtheorem{Corollary}[Theorem]{Corollary}
\newtheorem{Proposition}[Theorem]{Proposition}
\newtheorem{Lemma}[Theorem]{Lemma}
\theoremstyle{definition}
\newtheorem{Definition}[Theorem]{Definition}
\newtheorem{Note}[Theorem]{Note}
\newtheorem{Remark}[Theorem]{Remark}
\newtheorem{Example}[Theorem]{Example}
\newtheorem{Counterexample}[Theorem]{Counterexample}
\newtheorem{Assumption}[Theorem]{Assumption}
\parskip=.75pt
\makeatother

\title{Switching Diffusion Systems with Past-Dependent Switching and Countable State Space: Successful Couplings and Strong Ergodicity
}
\author{Fubao Xi\thanks{Frontier Interdisciplinary Domain, Beijing Institute of Technology, Zhuhai 519088, China, \tt
xifb@bit.edu.cn.} \and Yafei Zhai\thanks{School of Mathematics and Statistics, Beijing Institute of Technology, Beijing 100081, China, \tt yafeizhai@bit.edu.cn}
\and Chao Zhu\thanks{Department of Mathematical Sciences, University of Wisconsin-Milwaukee, Milwaukee, WI 53201, \tt zhu@uwm.edu.}}
 
\maketitle
\begin{abstract}
Motivated by applications in biology, ecology, and finance, this paper investigates a class of switching diffusion systems where the discrete component takes values in a countable state space with transition rates dependent on the history of the continuous component.  While history-dependent transition rates provide greater modeling flexibility for complex systems, they also introduce substantial analytical challenges due to the non-Markovian dynamics. In particular, successful coupling differs fundamentally from the delay-free setting: the meeting of two coupled processes does not immediately imply successful coupling, since the processes may separate again after coalescence as a result of their history-dependent switching. Consequently, successful coupling requires the two processes to remain identical for a sufficiently long period after their initial meeting. Under suitable conditions, we construct such a coupling, establish stability of the associated segment process in the total variation norm, and derive its strong ergodicity. Finally, we illustrate our theoretical results through an $N$-body mean-field model with past-dependent switching on a countable state space.

\bigskip

\bigskip

\noindent{\bf 2000 MR Subject Classification.} 60J25, 60J27, 60J60, 60J75.
\end{abstract}

\renewcommand{\baselinestretch}{1.18}

\section{Introduction}\label{sect-Intro}
Regime-switching stochastic differential equations have become a powerful framework for modeling and analyzing systems subject to random environmental changes, with applications in areas such as biological systems (\cite{BaoS-16,GLM-16,CCM-21,LLC-17}), control systems (\cite{WLXZ-23,Savku-24,Song-S-Z}), neural networks (\cite{Whitt-02, Huang-16}), mathematical finance and risk management (\cite{Zhang-01,Zhou-Yin,Elliott-Siu-10}), and population dynamics (\cite{Luo-Mao09,ZY-09c,SetLa-14}).
 These systems are characterized by abrupt transitions between different regimes, each governed by distinct dynamical laws. Due to their theoretical and practical importance, such  systems have attracted growing attention from researchers and practitioners in recent decades. We refer to \cite{MaoY} and \cite{YZ-10} for comprehensive investigations of such processes as well as their applications in areas such as mathematical biology, finance, and risk management. 

 The recent  decades have witnessed  extensive developments of the  theoretical framework of regime-switching   diffusions. To name just a few, 
\cite{XiYin-10} established asymptotic properties for nonlinear auto-regressive Markov processes with state-dependent switching, while \cite{XiYin-11} subsequently examined a class of jump-diffusions with state-dependent switching,
\cite{ZhuYB-15} derived Feynman-Kac formulas for regime-switching jump diffusions,  and  \cite{XiYin-15} obtained the exponential ergodicity for stochastic Li\'enard equations with state-dependent switching.  Further developments include the work of \cite{Shao-18}, which provided sufficient conditions for the existence and uniqueness of invariant measures for state-dependent regime-switching diffusion processes. \cite{NguyenY-20} investigated  almost sure and $L^p$ stability of stochastic functional differential equations with regime switching, while \cite{NNN-21} focused on their ergodicity and stability properties. \cite{TNY-22} analyzed stability in distribution of Markovian switching jump diffusions, and \cite{CWW-24} studied the weak convergence of stochastic functional diffusion systems with singularly perturbed regime switching.   We also refer to \cite{LMMY-20,LLLM-22} for investigations on delay feedback control strategies based on discrete-time observations of both the system and Markovian states.


Recent efforts have also focused on regime-switching (jump) diffusions where the switching component takes values in a countable state space. \cite{Shao-15} studied the existence of strong solutions and the Feller property for such systems. \cite{XiZ-17} established Feller and strong Feller properties and  exponential ergodicity for regime-switching jump diffusions with countable regimes. \cite{XiZW-21} explored stochastic damped Hamiltonian systems with state-dependent switching in a countably infinite regime space. In addition, \cite{Shao-17} investigated stabilization of regime-switching diffusion processes using feedback control based on discrete-time observations of both the state and switching processes, and \cite{ShaoX-19} extended these results to state-dependent regime-switching processes.

An important extension of regime-switching diffusions is to consider past-dependent switching, where the transition rates depend on the past trajectory of the process. This line of research  was initiated  in \cite{NguyenY-16},  which established existence and uniqueness results for a class of stochastic differential equations with past-dependent switching, and investigated their Markov and Feller properties. More recently, \cite{NNTY-23} examined the stability of stochastic functional differential equations with past-dependent random switching in a countably infinite state space. The past-dependent switching allows for more flexible modeling of systems where the transition rates are influenced by historical states, making it suitable for applications in ecology, biology,  finance, etc. For instance, in ecological systems, the transition rates of species populations may depend on their historical states, such as past population densities or environmental conditions. In financial markets, the transition rates of the general market trend (bull or bear) may depend on the  historical prices and/or the volatility patterns of the stocks. The past-dependent switching framework provides a more realistic representation of such systems, capturing the influence of historical states on current dynamics. 

Building upon the work of \cite{NguyenY-16}, we consider  the asymptotic properties for regime-switching  diffusion process $(X(t),\La(t))$ with past-dependent switching, in which the switching component $\La(t)$ has a countable state space. Since the switching rates of $\La(t)$ depends on the past states of $X$, the process $(X(t),\La(t))$ itself is not Markovian. However, the segment process $(X_t, \Lambda(t))$ is  Markovian as shown in \cite{NguyenY-16} (see also Lemma \ref{lem-Markov}). Consequently,  much of our analyses are focused on   $(X_t, \Lambda(t))$. In particular, we wish to establish strong ergodicity for $(X_t, \Lambda(t))$. Due to the complexity of the past-dependent switching, the analysis of such systems is significantly more challenging than that of standard regime-switching diffusions. The main difficulty arises from the fact that the past-dependent switching introduces inherently nonlinear memory effects, which invalidate the conventional approaches for (jump) diffusions with or without switching. Furthermore, the interaction between regime-switching and past-dependent dynamics complicates stability analysis, making it a non-trivial task to establish strong ergodicity. 

In this paper, we aim to address these challenges by developing a novel coupling methodology tailored for past-dependent switching diffusion processes. Our approach is inspired by the  coupling methods developed   in the study of (jump) diffusion processes such  as   \cite{Chen04}, \cite{ChenLi-89}, \cite{LindR-86}, \cite{PriolaW-06}, and \cite{Wang-10}. We construct a successful coupling for the past-dependent regime-switching diffusion processes, which allows us to derive strong ergodicity. This  advances previous work by \cite{Xi-13} on Markovian switching jump-diffusion processes and \cite{ShaoX-13,XiShao-13} on multidimensional diffusions with state-dependent switching. Specifically, we first introduce a coupling operator $\wdh \op$ in \eqref{e:wdh-op-reflection}, which is a combination of reflection and marching couplings when the discrete components are in the same regime and an independent coupling when the discrete components are  in different regimes; the operator $\wdh \op$ also incorporates the basic coupling for the discrete components. This gives rise to the coupling process $(X(t), \La(t), Y(t),  \Psi(t))$, which, in turn, leads to the segment process  $(X_t, \La(t), Y_t,  \Psi(t))$. 

 In the usual regime-switching diffusions without delay, the coupling time is typically defined as the first time when the two marginal processes  $(X(t), \La(t))$ and $(Y(t), \Psi(t))$ meet. For regime-switching diffusions with history-dependent switching, however, the situation is markedly different. Due to the delay, the processes $(X(t), \La(t))$ and $( Y(t),  \Psi(t))$ may separate again after their first meeting. Thus, coupling is successful only if the two processes remain identical for at least $r$ units of time after their initial coalescence, where $r$ denotes the length of the delay.  This, in turn, requires the switching components $\La$ and $\Psi$ to meet almost surely and, once they meet, to remain together for a sufficiently long period of time, even when starting from different initial conditions. To overcome the difficulties caused by the history-dependent switching, we exploit the temporal intervals during which the discrete components $\La$ and $ \Psi$ coincide. We then show that, almost surely, the processes $(X_t,\La(t))$ and $(Y_t, \Psi(t))$ coalesce during one of these intervals. This establishes the successful coupling of $(X_t,\La(t),Y_t, \Psi(t))$ asserted in Theorem \ref{success1}.

The successful coupling is then used in Corollary \ref{success3}  to show that the process $(X_t, \La(t))$ possesses a unique invariant measure $\pi$ and that the transition semigroup $P(t,(\phi,k),\cdot)$ converges to $\pi$  in the total variation norm as $t\to\infty$.   Furthermore, the coupling method allows us to show that the process $(X_t, \La(t))$ is strongly ergodic in the total variation norm. In other words, the exponential convergence of the transition semigroup $P(t,(\phi,k),\cdot)$ to the invariant measure $\pi$ is uniform with respect to the initial condition $(\phi, k)$. Moreover, an explicit exponential convergence rate is given. The details are spelled out in  Theorem \ref{thm:stroergo}. This  result  significantly broadens previous work on Markovian switching jump-diffusion processes and switching diffusions with state-dependent switching.

 Successful coupling and (strong) ergodicity for regime-switching (jump) diffusions have been extensively studied in the literature; see, e.g., \cite{ShaoX-13,Xi-13,XiShao-13,XiYin-11} and the references therein. To the best of our knowledge, the present work is the first to investigate strong ergodicity for switching diffusions with past-dependent switching through a successful coupling approach. Moreover, even in the delay-free setting, our conditions substantially extend existing results in the literature, including those in \cite{ShaoX-13,Xi-13,XiShao-13}. A detailed comparison with the existing literature is provided in Remark~\ref{Rem-comp-lit}.

The paper also establishes Feller property for the  segment process $(X_t, \Lambda(t))$ in Theorem \ref{thm-Feller}. The Feller property is a fundamental property that ensures  continuity of the transition semigroup and the strong Markov property of the underlying process. It also plays a key role in proving the ergodicity of the process. We show that the transition semigroup $P(t,(\phi,k),\cdot)$ associated with the past-dependent switching diffusion processes is continuous with respect to the initial condition $(\phi, k) $, thereby implying the Feller property. While the Feller property was established in \cite{NguyenY-16}, the proof was rather involved. In this paper,   we present a new and more streamlined approach by employing coupling methods to derive the Feller property.

As an application of the established results, we consider an $N$-body mean field model with past-state-dependent switching. This model  is a nontrivial  extension of the mean field model considered in \cite{XiYin-09} where the switching rates depend only on the current state of the system, not on the past. Additionally, instead of a finite state space for the switching component, we allow for a countable state space. Under some additional assumptions on the diffusion coefficients and switching rates, \cite{XiYin-09} established strong ergodicity for the underlying model. In this paper, we derive strong ergodicity for the more general model with less restrictions; the details are spelled out in Theorem \ref{thm:mf-strongergo}.  

The rest of the paper is arranged as follows. Section \ref{sect-Formulation} presents the precise formulation for past-dependent switching diffusion processes. The standing assumptions are also collected in Section \ref{sect-Formulation}. In Section \ref{sect-Feller}, we establish the Feller property for the corresponding regime-switching process $(X_t, \La(t))$. Sections \ref{sect-coupling} develops a successful coupling for the past-dependent switching diffusion processes. Section \ref{sect-stability} is devoted to the analysis of strong ergodicity  in
the sense of convergence in the total variation norm. Finally, as an application of the previous results, we discuss an $N$-body mean field model with past-state-dependent switching in Section \ref{sect-examples}.

\section{Formulation}\label{sect-Formulation}
To facilitate the presentation, we introduce some notation and definitions.   
For $x \in \rr^d$ and $\sigma
=(\sigma_{ij}) \in \rr^{d \times d}$, define
$$|x|=\bigl (\sum_{i=1}^{d} |x_{i}|^{2} \bigr )^{1/2}, \qquad
|\sigma|= \bigl ( \sum_{i,j=1}^{d} |\sigma_{ij}|^{2} \bigr
)^{1/2},$$ where $d$ is positive integer. Let  $\ss:=\{1,2,\cdots\}$ and denote by $d(\cdot, \cdot)$ the discrete metric on $\ss$. We next define a metric $\lambda (\cdot,\cdot)$ on $\rr^{d}
\times \ss$ by $$\lambda \bigl((x,m), (y,n) \bigr)=|x-y|+d(m,n), \quad \forall \, (x,m), (y,n) \in \rr^{d} \times \ss.$$
 Let ${\cal
B}(\rr^{d} \times \ss)$ be the Borel $\sigma$-algebra on $\rr^{d}
\times \ss$. Then $(\rr^{d} \times \ss, \lambda (\cdot,\cdot), {\cal
B}(\rr^{d} \times \ss))$ is a locally compact and separable metric
space. Next, denote by $\C:=\C([-r, 0], \rr^d)$ the set of $\rr^d$-valued continuous functions, where $r$ is a fixed positive number. Equip $\C$ with the super norm, i.e.,  $\|\phi\|:=\sup\{|\phi(t)|: t\in [-r, 0]\}$ for any $\phi \in \C$. For any $y\in \C$ and  $t\ge 0$, denote by $y_t$ the so-called segment function (or memory segment function) $y_t:=\{y(t+s): -r\le s\le 0\}$. {Let $E:=\mathcal C\times\mathbb S$ and denote by   $\mathcal{P}(E)$   the collection of all probability measures on $(E, \mathcal{B}(E))$. The  space of bounded and continuous functions on $E$ is denoted by $C_b(E)$. Throughout the paper, $\1_A$ denotes the indicator function of a set $A$ while $I$ is the identity matrix of appropriate dimension.} Finally,  $(\Omega,  {\cal F}, \{{\cal
F}_{t}\}_{t\ge 0},  \bP)$ is a complete probability space with a
filtration $\{{\cal F}_{t}\}_{t\ge 0}$ satisfying the usual
conditions (i.e., it is increasing and right continuous and ${\cal
F}_{0}$ contains all $\bP$-null sets).

To formulate our model, let $(X(t),\Lambda(t))$ be a switching diffusion process with past-dependent switching on $\rr^d \times \ss$. The
first component $X(t)$ satisfies the following stochastic
differential equation (SDE)
\begin{equation}\label{eq:X}
\d X(t)=
\displaystyle b(X(t),\Lambda(t))\d t + \sigma (X(t),\Lambda(t))\d B(t),
\end{equation}
where 
$b(x,k): \R^d\times \ss\mapsto \R^d$ and $\sigma(x,k): \R^d\times \ss\mapsto \R^{d\times d}$ are Borel measurable functions, and 
 $B(\cdot)$ is a $ {d}$-dimensional standard Brownian motion.
The second component $\Lambda(t)$ is a discrete random 
process with  
state space $\ss$ such that
\begin{equation}\label{eq:Lambda}
\P \{\Lambda(t+\Delta)=l | \Lambda(t)=k, X_t=\phi\} =
\begin{cases}q_{kl}(\phi) \Delta +o(\Delta),
\ad \, \, \hbox{if}\, \, l \ne k, \\
1+q_{kk}(\phi) \Delta +o(\Delta), \ad \, \, \hbox{if}\, \, l = k,
\end{cases} 
\end{equation} provided $\Delta \downarrow 0$, where the switching tensity matrix
$Q(\phi)= \bigl(q_{kl}(\phi)\bigr)$ depends on the past trajectory of $X(t)$ in the interval $[t-r, t]$, that is, the  
 segment process $X_t$. As
usual, we assume that $\sum_{l\in \ss}q_{kl}(\phi)=0$ for all $\phi\in \C$ and $k\in\ss$. 

For the existence and uniqueness of a strong solution
$(X(t),\Lambda(t))$ to  the system (\ref{eq:X})--(\ref{eq:Lambda}), we make the following assumption:

\begin{Assumption}\label{I1} \begin{itemize}
\item[(i)]  For each $k\in\ss$, the functions $b(\cdot,k)$ and $\sg(\cdot,k)$ satisfy the local Lipschitz condition. That is, for every $R>0$, there exists a positive constant $L_{R}$ such that 
\begin{equation}
\label{e:b/sig-Lip}
|b(x, k) - b(y, k) | + |\sigma(x, k)-\sigma(y, k)|\le L_{R} |x-y|, 
\end{equation} for all $x,y\in \R^{d}$ with $|x|\vee|y| \le R$   and $ k\in \ss$. 
 \item[(ii)] {For each $k\in \ss$, there exist a positive constant $\gamma_k$ and a nonnegative function $V_k\in C^2(\R^d)$ such that   \begin{align*} & \lim_{R\to\infty}  \inf\{V_k(x): x\in \R^d, |x|\ge R\}  = \infty, \\
  \LL_k V_k(x) & := \frac{1}{2} \hbox{tr}\bigl(\sigma(x,k)\sigma^\top(x,k)D^2V_k(x)\bigr) + \lan b(x,k), D V_k(x)\ran \le \gamma_k (1+V_k(x)), 
 \end{align*}  for all $x\in \R^d$,   where $D V_k(x)$ and $D^2 V_k(x)$ denote the gradient and Hessian matrix of $V_k$ with respect to $x$, respectively.}
\item[(iii)] \ There exists a constant $H>0$ such that
\begin{equation}\label{(1.6)}
\sup\{q_k(\phi):\phi\in\C, k\in \ss \} \le H<\infty, 
\end{equation}\end{itemize}
\end{Assumption}

Under Assumption~\ref{I1}, for each $k\in \ss$,   the   SDE
\begin{equation}\label{eq:Xk}
\d X^{(k)}(t)=b(X^{(k)}(t),k)\d t+\sigma(X^{(k)}(t),k)\d B(t),\quad  X^{(k)}(0) =x \in \R^d \end{equation}
admits a unique non-explosive strong solution $X^{(k)}(\cdot)$. Furthermore,   using the interlacing procedure together with exponential killing method as in \cite{XiZ-17,NguyenY-16}, or \cite{APPLEBAUM}, we can prove that for given 
initial data $(\phi,k)\in \C\times \ss$, there exists a unique non-explosive solution such that $X(t)=\phi(t)$ for $t\in [-r,0]$, $\Lambda(0)=k$, and $(X(t),\Lambda(t))$ satisfies system (\ref{eq:X})--(\ref{eq:Lambda}). For the simplicity of presentation, we will assume that Assumption~\ref{I1} holds in the rest of the paper.
In general, the solution $(X(t),\Lambda(t))$ to system (\ref{eq:X})--(\ref{eq:Lambda}) is not a Markov process, but $(X_t,\Lambda(t))$ is a strong Markov  process under Assumption \ref{I1}; see Theorem 4.1 of \cite{NguyenY-16}. 

 We finish this section with two examples. The first example is taken from \cite{NguyenY-16} and concerns with the dynamics of interacting micro- and macro-species.  Denote by $X(t)\in [0,\infty)$ and $\La(t)\in \ss$ the bio-density and the population size of the micro- and macro-species at time $t$, respectively. The life cycle of a micro species is usually very short, whereas the reproduction  of a macro species is  non-instantaneous and typically depends on the period of time from egg formation to hatching, say, $r$. Therefore, it is natural to assume that the transition rates of $\La(t)$ depend on the  history of the processes. This gives rise to the model \eqref{eq:X}--\eqref{eq:Lambda}, where  \eqref{eq:X} is specified as a stochastic logistic model: $$\d X(t)=\disp X(t)\bigl(a(\La(t))-b(\La(t))X(t)\bigr)\d t+\sigma (\Lambda(t))X(t)\d B(t),$$ where $a(k)$, $b(k)$ and $\sg(k)$ are positive constants for each $k\in \ss$.

The second example considers a regime-switching Black-Scholes model featuring  past-dependent  transition rates. This is motivated by the empirical studies in \cite{LiuMW:12,SchavN:97}, which show that market trends often depend on historical stock prices. Suppose the price of a risky asset $X(t)$ follows the dynamics
$$\d X(t)=\disp X(t)\bigl(\mu(\La(t))\d t+\sigma (\Lambda(t))\d B(t)\bigr),$$ where the switching component $\La(t) \in \{1, 2\}$ represents the market trend, which can be either bullish (if $\La(t)=1$) or bearish (if $\La(t)=2$), and for $k=1, 2$, $\mu(k)$ and $\sg(k)$ are positive constants. Unlike the traditional regime-switching Black-Scholes models considered in \cite{Elliott-Siu-10,WLXZ-23,Zhang-01,Zhou-Yin}, we allow the transition rates of $\La(t)$ to depend on the past trajectory of $X(t)$ as  in \eqref{eq:Lambda}, thereby capturing the influence of historical asset prices on market regimes and allowing for more realistic modeling of financial markets. For example, defining  $f(\phi) : = \log \frac{\|\phi\|}{\phi(0)}$ for any $\phi \in \mathcal C([-r, 0]; \R_+)$, in which $\|\phi\| : = \sup_{s\in [-r,0]} |\phi(s)|$,   
 we may  specify the transition rates in \eqref{eq:Lambda} as  \begin{align*}
  q_{12}(\phi) : = q_{1} + \frac{f(\phi)}{1+ f(\phi)}, \qquad q_{21}(\phi) : = q_2 \bigg(1- \frac{f(\phi)}{1+ f(\phi)}\bigg),
 \end{align*} where $q_1, q_2$ are positive constants.  Under such a specification, if the current  price $X(t)$ is much lower than its running maximum $\sup_{s\in[0 \vee(t-r), t]} X(s)$, then $f(X_t)$ and thus the transition rate $q_{12}(X_t)$ increase; indicating a more likely transition from bull to bear markets. Likewise, if the current price is close to its running maximum, a transition from a bear to bull market is more likely. 
 
 Alternatively, we may develop a regime-switching Black-Scholes model where the transition rates depends on the moving averages of the underlying risky asset.

\section{Feller Property}\label{sect-Feller}
This section is devoted to establishing the Feller property for the process   $(X_t,\La(t))$. To this end, we impose the following condition: 
\begin{Assumption}\label{assump-Q}
There exists a   concave function $\gamma: [0,\infty)\mapsto [0, \infty)$ with  $\lim_{x\downarrow 0}\gamma(x) = \gamma(0) = 0$ 
such  that for every $k\in \ss$ and  $R >0$, 
\begin{equation}
\label{e:Q-cont}
\sum_{j\in \ss\setminus\{k\}} |q_{kj}(\phi) - q_{kj}(\psi) | \le  { \kappa_{R} \gamma(  \norm{\phi-\psi})}, \quad \forall \phi, \psi \in \C \text{ with }\|\phi\|\vee \|\psi\| \le R,
\end{equation} where $\kappa_{R}$ is a positive constant. 
\end{Assumption}

\begin{Theorem}\label{thm-Feller}
Suppose Assumptions \ref{I1} and \ref{assump-Q} hold. Then for any $f\in C_{b}(\C\times \ss; \R)$ and $t> 0$, the function \begin{equation}\label{e:Ptf-defn}
P_{t} f(\phi, k): =\E[f(X^{\phi, k}_{t}, \La^{\phi, k}(t))], \quad (\phi, k)\in \C\times \ss,
\end{equation} is continuous. Moreover, we have $\lim_{t\downarrow 0}P_{t} f(\phi, k) =f(\phi, k) $.
\end{Theorem}
We will prove Theorem \ref{thm-Feller} by  the coupling method. To this end,  instead of the infinitesimal equation \eqref{eq:Lambda}, we need to describe the evolution of $\Lambda(t)$ using a stochastic differential equation. We first construct a family $\{\Delta_{ij}(\phi), j\not=i\in \ss, \phi\in \C\}$ of left-closed, right-open intervals on the positive half real line. In addition, for each $j\neq i$, the interval $\Delta_{ij}(\phi)$ has length $q_{ij}(\phi)$.  We set
$\Delta_{ij}(\phi) = \emptyset$ if $q_{ij}(\phi) = 0$, $i \not= j$. We refer to \cite{NguyenY-16} for more details on the construction of these intervals.  Next we define a function
$h: \C \times \ss \times \R_{+} \mapsto \R$ by
\begin{equation}\label{eq:h-def} h(\phi,i,z)=\sum_{j\in\ss\setminus\{i\}} (j-i)
\1_{\{z\in\Delta_{ij}(\phi)\}}.
\end{equation}
That is,  
$h(\phi,i,z)=j-i$ if $z\in\Delta_{ij}(\phi)$; otherwise $h(\phi,i,z)=0$. 
 We can now consider the following stochastic differential equation for $\La(t)$:
\begin{equation}\label{2-eq:ju} \d \La(t) = \int_{\R_{+}} h(X_{t-},\La(t-),z) {\mathfrak p}(\d t,\d z),\end{equation}
where ${\mathfrak p}(\d t,\d z)$ is a Poisson random measure on $[0,\infty) \times \R_{+}$ with intensity $\d t
\times \m(\d z)$, and $\m(\cdot)$ is the Lebesgue measure on $\R_+$. The Poisson random
measure ${\mathfrak p}(\cdot, \cdot)$ is independent of the Brownian motion
$B(\cdot)$. 

\begin{Lemma}\label{lem-La_SDE}{\em \cite[Lemma A.1]{NguyenY-16}}
  Suppose Assumptions \ref{I1} and \ref{assump-Q} hold, then the system \eqref{eq:X}--\eqref{eq:Lambda} is equivalent to the system \eqref{eq:X}--\eqref{2-eq:ju}. 
\end{Lemma}

 In view of Lemma \ref{lem-La_SDE} and Theorem 3.5 of \cite{NguyenY-16}, the system \eqref{eq:X}--\eqref{2-eq:ju} is well-posed in the strong sense. Consequently, we can assume that the filtration $\{{\cal F}_{t}\}_{t\ge 0}$ is generated by the initial condition $(X_{0},\Lambda(0))$, the Brownian motion $B(\cdot)$, and the Poisson random measure ${\mathfrak p}(\cdot, \cdot)$,  augmented by the $\P$-null sets. In other words, $\{\F_{t}\}$ is generated as follows. First, we let \begin{displaymath}
\wdt{\mathcal F}_{t}: = \sigma(X_{0},\La({0}), B(s), \mathfrak p((0, s]\times U): 0\le s \le t, U\in \B(\R_{+})), \ 0\le  t< \infty,
\end{displaymath} and $\wdt{\mathcal F}_{\infty}: = \sigma\big(\bigcup_{t\ge 0} \wdt{\mathcal F}_{t}\big) $, as well as the collection of null sets:
\begin{displaymath}
\mathcal N:= \{N \subset \Omega: \text{there exists a }G\in \wdt{\mathcal F}_{\infty} \text{ with } N \subset G \text{ and }\P(G) = 0\}.
\end{displaymath} We then create the augmented filtration
\begin{displaymath}
\F_{t} : = \sigma(\wdt{\mathcal F}_{t} \cup \mathcal N), \ 0\le t < \infty; \ \ \F_{\infty} : = \sigma\bigg(\bigcup_{t\ge 0} \F_{t}\bigg).
\end{displaymath} In particular, the process $(X(t),\Lambda(t))$ (and hence the segment process $(X_t,\Lambda(t))$ as well)  is adapted to the filtration $\{{\cal F}_{t}\}_{t\ge 0}$.  
Furthermore, the following lemma shows that the   process $(X_t,\La(t))$ is   Markovian with respect to the filtration $\{\F_t\}_{t\ge 0}$. 

 \begin{Lemma}\label{lem-Markov}
		Suppose Assumptions \ref{I1} and \ref{assump-Q} hold, and denote by $(X(t),\La(t))$ the solution to \eqref{eq:X}--\eqref{2-eq:ju} with initial condition $(\phi, k) \in \C\times\ss$. Then $(X_t,\Lambda(t))$ is a Markov  process. In other words, for any $t\ge s\geq 0$ and   any bounded Borel measurable function $\bar h:E\to \mathbb R$,  we have
    $$\E[ \bar h(X_t, \Lambda(t))|\mathcal F_s] =\E [ \bar h(X_t, \Lambda(t))|(X_s,\Lambda(s))], \quad\text{a.s.}$$   
\end{Lemma}
\para{Proof.} 
For any initial time $s\geq0$, we consider the following system for $t\geq s$,
$$X^{s,\phi,k}(t)=\phi(0)+\int_{s}^{t} b(X^{s,\phi,k}(u),\Lambda^{s,\phi,k}(u))\d u +\int_{s}^{t} \sigma (X^{s,\phi,k}(u),\Lambda^{s,\phi,k}(u))\d B(u),$$
and
$$\La^{s,\phi,k}(t) =k+\int_{s}^{t} \int_{\R_{+}} h(X^{s,\phi,k}_{u-},\La^{s,\phi,k}(u-),z) {\mathfrak p}(\d u,\d z),$$
with the initial condition
$$X^{s,\phi,k}(m)=\phi(m-s), \text{for $m\in[s-r, s]$ and $\phi\in\mathcal{C}$}.$$
Under Assumptions \ref{I1} and \ref{assump-Q}, the system admits a unique non-explosive strong solution $(X^{s,\phi,k}(t), \Lambda^{s,\phi,k}(t))$ and we denote the corresponding segment process by $(X^{s,\phi,k}_t, \Lambda^{s,\phi,k}(t))$.

Let $\mathcal G_s=\sigma\{B(u)-B(s), {\mathfrak p}(u,U)-{\mathfrak p}(s,U):u\geq s, U\in \mathcal B(\mathbb R_+)\}$. Clearly, $\mathcal G_s$ is independent of $\mathcal F_s$. Moreover, the process $(X^{s,\phi,k}_t, \Lambda^{s,\phi,k}(t))$ depends completely on the increments $B(u)-B(s), {\mathfrak p}(u,\d z)-{\mathfrak p}(s,\d z)$ for $u\geq s$ and so is $\mathcal G_s$-measurable. Hence $(X^{s,\phi,k}_t, \Lambda^{s,\phi,k}(t))$ is independent of $\mathcal F_s$ for all $t\geq s$. On the other hand, the strong uniqueness of the solution implies that 
$$X(u)=X^{s,X_s,\La(s)}(u)\ \text{ for any $u\geq s-r$},$$
and
$$\La(u)=\La^{s, X_s,\La(s)}(u)\ \text{ for any $u\geq s$}.$$
Thus, with probability one, we have $X_t=X^{s,X_s,\La(s)}_t$ for any $t\geq s$.

For any bounded Borel measurable function $\bar h:E\to \mathbb R$, by Proposition 1.12 of \cite{DaPratoZ-14}, we have
\begin{align*}
\E[ \bar h(X_t, \La(t))|\mathcal F_s]&=\E[\bar h(X^{s,X_s,\La(s)}_t, \La^{s, X_s,\La(s)}(t))|\mathcal F_s]\\
	&=\E[ \bar h(X^{s,\phi,k}_t, \La^{s, \phi,k}(t))|(\phi,k)=(X_s,\La(s))]\\
	&=\E[ \bar h(X_t, \La(t))|(X_s,\La(s))].
\end{align*}
This establishes  the Markov property.
\qed

In view of Lemma \ref{lem-La_SDE}, we now consider the basic coupling for the process $(X,\La)$: 
\begin{equation}\label{e-basic-coupling}
\begin{cases}
\d X(t)=b(X(t), \La(t))dt+\sigma(X(t),\La(t))\d B(t), \\
\d\La(t) = \int_{\R_{+}} h(X_{t-}, \La(t-),z) \mathfrak p(\d t,\d z),\\
\d\wdt  X(t)=b(\wdt X(t), \wdt \La(t))dt+\sigma(\wdt  X(t),\wdt \La(t))\d B(t),\\
\d\wdt \La(t) = \int_{\R_{+}} h(\wdt X_{t-}, \wdt \La(t-),z) \mathfrak p(\d t,\d z),
\end{cases}
\end{equation} with initial condition $ (\phi,k, \psi, k)\in \C\times \ss\times \C\times \ss$.  One can show that a unique solution of \eqref{e-basic-coupling} exists. For $f \in C_{c}^{2} (\R^{d} \times \M\times \R^{d} \times \M;\R)$, 
we define the operator associated with \eqref{e-basic-coupling} as
\begin{equation}
\label{e:op-basic-coupling}
\wdt \op f(\phi, i , \psi, j) : =    \wdt \varOmega_{\text{d}}   f(\phi(0),i, \psi(0),j) +   \wdt \varOmega_{\text{s}}f(\phi,i, \psi,j), \quad \forall (\phi, i , \psi, j) \in \C\times \ss\times \C\times \ss,
\end{equation} where $ \wdt \varOmega_{\text{d}}$   and    $ \wdt \varOmega_{\text{s}}$ are defined as follows.
For $x,z\in \R^{d}  $ and $i,j\in \M $, we set 
\begin{align}\label{eq:double-a}a(x,i,z,j) & =\begin{pmatrix}
  \sigma(x, i)\\ \sigma(z,j)\end{pmatrix} \begin{pmatrix}
    \sigma^\top(x, i) & \sigma^\top(z,j)
  \end{pmatrix}= \begin{pmatrix}
a(x,i) & \sigma (x,i) \sigma^\top (z,j) \\
\sigma (z,j) \sigma^\top (x,i) & a(z,j)
\end{pmatrix},\end{align} and \begin{align} \label{eq:double-b}
b(x,i,z,j) =\begin{pmatrix}
b(x,i)\\
b(z,j) \end{pmatrix}, \end{align} where $a(x,i)= \sigma(x,i)\sigma^\top(x,i)$. Then we define for any $(x,i,z,j)\in \R^d\times \ss\times \R^d\times \ss$,
\begin{align}\label{eq-Omega-d-defn}
& \wdt {\varOmega}_{\text{d}}f(x,i,z,j): =\frac
{1}{2}\hbox{tr}\bigl(a(x,i,z,j)D^{2}f(x,i, z,j)\bigr)  +\langle
b(x,i,z,j), D f(x,i, z,j)\rangle,
\end{align} where $D f(x,i, z,j)= (D_{x } f(x,i,z,j), D_{z} f(x,i,z,j))$  and $D^{2}f(x,i,z,j)$ denote the gradient and   the Hessian matrix of   $f$ with respect to the  variables $x$ and $z$, respectively. Finally, the operator  $\wdt \varOmega_{\text{s}}  f$ is defined as follows:  
\begin{align}
\label{eq-Q(x)-coupling}
\nonumber \wdt \varOmega_{\text{s}}  f(\phi,i,\psi,j) & : =
 \sum_{l\in\M}[q_{il}(\phi)-q_{jl}(\psi)]^+(   f(\phi(0) ,l, \psi (0), j)-  f(\phi(0),i, \psi(0), j))\\
\nonumber &\quad\ +\sum_{l\in\M}[q_{jl}(\psi) -q_{il}(\phi)]^+( f(\phi(0),i, \psi(0), l)- f(\phi(0),i, \psi(0), j) )\\
&\quad \  +\sum_{l\in\M}[ q_{il}(\phi) \wedge q_{jl}(\psi) ](  f(\phi(0),l, \psi(0), l)-f(\phi(0),i, \psi(0), j)).
 \end{align}
 
We are now ready to prove Theorem \ref{thm-Feller}.
\begin{proof}[Proof of Theorem \ref{thm-Feller}] The assertion that $\lim_{t\downarrow 0}P_{t} f(\phi, k) =f(\phi, k) $ is obvious because the process $(X, \La)$ is c\`adl\`ag and $f$ is bounded and continuous. To establish the continuity of the function $P_{t}f$ of \eqref{e:Ptf-defn}, we consider the coupling process $(\wdt X, \wdt \La, X, \La)$ given by \eqref{e-basic-coupling}  with initial condition $ (\phi,k, \psi, k)\in \C\times \ss\times \C\times \ss$; here we take $\wdt\La(0) = \La(0) $ because $\ss$ has a discrete topology. 
For each $R> 0$, let $$\tau_{R}: = \inf\{t\ge 0: \|X_{t}\| \vee  \|\wdt X_{t}\| \ge R\}.  $$ Also denote by  $$ \zeta: = \inf\{t\ge 0: \La(t) \neq \wdt\La(t)\}$$   the first time when the discrete components $\La$ and $\wdt \La$ differ. We have $\P(\zeta > 0) =1$. For any $s\in [0, \zeta)$, we have 
\begin{align*}
\wdt X(s) - X(s)&  = \psi(0) - \phi(0) + \int_{0}^{s} [b(\wdt X(r), \La(r))-b(  X(r),   \La(r))] \d r \\  & \quad +  \int_{0}^{s} [\sigma(\wdt X(r), \La(r))-\sigma( X(r),   \La(r))] \d B(r).
\end{align*}  Then it follows from the local Lipschitz condition for $b(\cdot, k)$ and $\sigma(\cdot, k), k\in \ss$ that for any $t\in [ 0, T]$ with $T>0$ being temporarily fixed, 
\begin{align*} 
\E& \bigg[\sup_{0\le s \le t\wedge \tau_{R}\wedge \zeta} |\wdt X(s) - X(s)|^{2}\bigg] \\& \le 3 |\psi(0) - \phi(0)|^{2} + 3 t \E\bigg[\int_{0}^{t\wedge \tau_{R}\wedge \zeta} |b(\wdt X(r), \La(r))-b(  X(r),   \La(r))|^{2}\d r \bigg] \\
 & \quad + 12 \E\bigg[\int_{0}^{t\wedge \tau_{R}\wedge \zeta} |\sigma(\wdt X(r), \La(r))- \sigma(  X(r),   \La(r))|^{2}\d r \bigg] \\
 & \le 3 |\psi(0) - \phi(0)|^{2}  + 3(t+4) L_{R}^{2} \E\bigg[\int_{0}^{t\wedge \tau_{R}\wedge \zeta} |\wdt X(r) - X(r)|^{2} \d r\bigg]\\
 & \le 3 |\psi(0) - \phi(0)|^{2}  + 3(T +4) L_{R}^{2} \int_{0}^{t} \E\bigg[\sup_{ 0\le u \le r\wedge \tau_{R}\wedge \zeta }  |\wdt X(u) - X(u)|^{2} \bigg]\d r. 
\end{align*} The Gronwall inequality then implies that 
\begin{displaymath}
\E  \bigg[\sup_{0\le s \le t\wedge \tau_{R}\wedge \zeta} |\wdt X(s) - X(s)|^{2}\bigg] \le 3 |\psi(0) - \phi(0)|^{2} e^{3(T +4) L_{R}^{2} t}.
\end{displaymath}
Then it follows that for any $s\in[0, t]$, \begin{align} \label{eq:tX-X=0}
 \E[\| \wdt X_{s \wedge \tau_{R}  \wedge\zeta}- X_{s \wedge \tau_{R}  \wedge\zeta} \|] & \le \|\phi-\psi\| \vee \E \bigg[\sup_{0\le s \le t\wedge \tau_{R}\wedge \zeta} |\wdt X(s) - X(s)|\bigg]
 \le C_{T,R}\|\phi-\psi\|, 
\end{align} where $C_{T, R} > 0$ is a constant depending on $T$ and $R$.

Consider the function $\Xi(\phi(0),k,\psi(0),l): = \1_{\{k\neq l\}}$ for $(\phi, k, \psi, l) \in \C\times \ss\times \C\times \ss$. It follows directly  from the definition of $\wdt \op$ that $$\wdt   \op \Xi(\phi,k,\psi,l) = \wdt \varOmega_{\text{s}} \Xi(\phi,k,\psi,l) \le 0, \text{ if } k \neq l.$$   When $k =l$, we have  from \eqref{e:Q-cont} that
\begin{align*}
   \wdt   \op \Xi(\phi,k,\psi,l) & =\wdt \Omega_{\text{s}} \Xi( \phi,k,\psi,k) & \\
   &   = \sum_{ i\in \M} [q_{ki}(\phi)-q_{ki}(\psi)]^+( \1_{\{i \neq k \}} - \1_{\{ k\neq k\}}) 
   \\\nonumber  &\qquad
 +\sum_{i \in\M}[q_{ki}(\psi)  -q_{ki}(\phi)]^+(  \1_{\{i \neq k \}} - \1_{\{ k\neq k\}} ) \\
 & \le \sum_{i \in\ss\setminus\{ k\}} \abs{q_{ki}(\phi)-q_{ki}(\psi)}
 \le  \kappa_{R} \gamma( \norm{\phi-\psi}),
\end{align*} for all $\phi,\psi\in \C$ with $\|\phi\|\vee\|\psi\| \le R$.
Hence a combination of the previous two displayed equations leads to  \begin{equation}
\label{eq-coupling est 1}
\wdt \op \Xi(\phi,k,\psi,l) \le   { \kappa_{R} \gamma(  \norm{\phi-\psi})}
\end{equation}
for all  $(\phi, k, \psi, l) \in \C\times \ss\times \C\times \ss$ with  $\|\phi\| \vee \|\psi\| \le R$.

 Note that $\zeta \leq t \wedge \tau_{R}  $ if and only if  $\wdt {\La}(t \wedge \tau_{R}  \wedge \zeta) \neq  \La(t \wedge \tau_{R}  \wedge \zeta)$. In addition, $\La(s) = \wdt\La(s)$ for all $s\in [0, t \wedge \tau_{R}  \wedge\zeta) $. Thus  we can use \eqref{eq-coupling est 1} to  compute
  \begin{align*}
	\mathbb{P}  \{\zeta \leq t \wedge \tau_{R}   \} 
	& = \E[\Xi(\wdt  X(t \wedge \tau_{R}  \wedge \zeta), \wdt  \La(t \wedge \tau_{R}  \wedge \zeta),    X(t \wedge \tau_{R}  \wedge \zeta),\La(t \wedge \tau_{R}  \wedge \zeta))]\\
		&=\Xi(\wdt  x,k,x,k) + \mathbb{E}\bigg[\int_{0}^{t \wedge \tau_{R}  \wedge\zeta}\wdt \op \Xi(\wdt  X_{s},\wdt \La(s), X_{s},\La(s)) \d s\bigg]\\
	&\leq \kappa_{R}\mathbb{E}\bigg[\int_{0}^{t\wedge \tau_{R}  \wedge\zeta }  { \gamma (\|\wdt X_{s}- X_{s}\|)  } \d s\bigg]  \\ &  
	\leq \kappa_{R}\int_{0}^{t} \mathbb{E}\big[ \gamma(\|\wdt X_{s \wedge \tau_{R}  \wedge\zeta}- X_{s \wedge \tau_{R}  \wedge\zeta}\| ) \big ] \d s   \\ &
	\le   \kappa_{R}\int_{0}^{t}{  \gamma\big( \mathbb{E}[\|\wdt X_{s \wedge \tau_{R}  \wedge\zeta}- X_{s \wedge \tau_{R}  \wedge\zeta}\|]\big)  } \d s,
\end{align*}
where we used the concavity of $\gamma$ to derive  the last inequality. 
Note that  $\|\wdt X_{s \wedge \tau_{R}  \wedge\zeta}- X_{s \wedge \tau_{R}  \wedge\zeta}\| \le 2 R$ for all $s\in [0, t]$. Then it follows from 
\eqref{eq:tX-X=0} and the bounded convergence theorem that 
\begin{equation}\label{e:p1=0}
\lim_{\|\phi   -\psi\| \to 0}\P\{\zeta\le t \wedge \tau_{R}   \}=0.
\end{equation} 
Now for any $f\in  C_{b}(\C\times \ss; \R)$ and $t\ge 0$, we compute
\begin{align} \label{eq-feller}
\nonumber |&P_{t} f(\psi, k) - P_{t} f(\phi, k)|\\\nonumber & = |\E[f(\wdt X_{t}, \wdt\La(t)) - f(X_{t}, \La(t))]| \\
\nonumber  & \le \E[|f(\wdt X_{t}, \wdt\La(t)) - f(X_{t}, \La(t))| \1_{\{\tau_{R} \le t \}}] + \E[|f(\wdt X_{t}, \wdt\La(t)) - f(X_{t}, \La(t))| \1_{\{\tau_{R} > t \}}] \\
\nonumber  & \le 2 \|f\|_{\infty} \P\{\tau_{R} \le t \} +  \E[|f(\wdt X_{t\wedge \tau_{R}}, \wdt\La(t\wedge \tau_{R})) - f(X_{t\wedge \tau_{R}}, \La(t\wedge \tau_{R}))| \1_{\{\tau_{R} > t \}}] \\
\nonumber  & \le 2 \|f\|_{\infty} \P\{\tau_{R} \le t \} + \E[|f(\wdt X_{t\wedge \tau_{R}}, \wdt\La(t\wedge \tau_{R})) - f(X_{t\wedge \tau_{R}}, \La(t\wedge \tau_{R}))| \1_{\{\zeta \le \tau_{R}  \wedge t \}}] \\
\nonumber  & \quad + \E[|f(\wdt X_{t\wedge \tau_{R}}, \wdt\La(t\wedge \tau_{R})) - f(X_{t\wedge \tau_{R}}, \La(t\wedge \tau_{R}))| \1_{\{\zeta >  \tau_{R}  \wedge t \}}] \\
 & \le  2 \|f\|_{\infty} \P\{\tau_{R} \le t \} +2 \|f\|_{\infty}\P\{\zeta\le t \wedge \tau_{R}   \} \\
\nonumber  & \quad + \E[|f(\wdt X_{t\wedge \tau_{R}\wedge\zeta}, \La(t\wedge \tau_{R}\wedge\zeta)) - f(X_{t\wedge \tau_{R}\wedge\zeta}, \La(t\wedge \tau_{R}\wedge\zeta))| \1_{\{\zeta >  \tau_{R}  \wedge t \}}].
\end{align} For any $\e >0$, since $\tau_{R} \to \infty$ a.s. as $R \to \infty$, we can choose an $R$ sufficiently large so that $ 2 \|f\|_{\infty} \P\{\tau_{R} \le t \} < \frac\e3$. Likewise, \eqref{e:p1=0} allows us to choose $\|\psi-\phi\|$ sufficiently small so that $2 \|f\|_{\infty}\P\{\zeta\le t \wedge \tau_{R}   \} < \frac\e3$. In a similar manner, we can use \eqref{eq:tX-X=0}, the continuity of $f$, and the bounded convergence theorem to conclude  that for all $\|\psi-\phi\|$ sufficiently small, we have\begin{displaymath}
\E[|f(\wdt X_{t\wedge \tau_{R}\wedge\zeta}, \La(t\wedge \tau_{R}\wedge\zeta)) - f(X_{t\wedge \tau_{R}\wedge\zeta}, \La(t\wedge \tau_{R}\wedge\zeta))| \1_{\{\zeta >  \tau_{R}  \wedge t \}}]  < \frac\e3.
\end{displaymath} Plugging these observations into \eqref{eq-feller}, we derive \begin{displaymath}
|P_{t} f(\psi, k) - P_{t} f(\phi, k)| < \e.
\end{displaymath} Since $\e> 0$ is arbitrary, we obtain the desired continuity of the function $P_{t}f$ defined in \eqref{e:Ptf-defn}. This completes the  proof. 
\end{proof}

\begin{Remark}
  While Feller property was established in \cite{NguyenY-16}, we have provided a new proof based on the coupling method. The advantage of our proof is that it is more direct and does not require the use of the change of measure approach.
\end{Remark}

\begin{Remark}
 Since the trajectories of $(X_t,\La(t))$ are c\`adl\`ag,  it follows from Theorem \ref{thm-Feller} and Lemma \ref{lem-Markov} that $(X_t, \Lambda(t))$ is a strong Markov process.
\end{Remark}
  
\section{Successful Coupling}\label{sect-coupling}
%

In order to construct a coupling process $(X(t),\La(t),Y(t), \Psi(t))$ for two copies $(X(t),\La(t))$ and $(Y(t), \Psi(t))$ of   solutions to the  system \eqref{eq:X}--\eqref{eq:Lambda}. For $(x,k,y,l) \in \rr^d \times \ss \times \rr^d \times \ss$, set a  $2d \times 2d$ matrix as follows:
\begin{align}\label{sgxkyl} \nonumber  \tau(x,k,y,l)& =  {\mathbf{1}}_\triangle (k,l){\mathbf{1}}_{\{x\ne y\}}
\begin{pmatrix}
  \sg(x,k) & 0\\
  \sg(y,k)H(x,y,k) & 0
\end{pmatrix}  
 \,+{\mathbf{1}}_\triangle (k,l){\mathbf{1}}_{\{x=y\}}
\begin{pmatrix}
  \sg(x,k) & 0\\
  \sg(x,k) & 0
\end{pmatrix} \\
\ad \ \ \quad +{\mathbf{1}}_{\triangle^c}
(k,l)\begin{pmatrix}
  \sg(x,k) & 0 \\
  0 & \sg(y,l)
\end{pmatrix},\end{align}
where $\triangle=\{(k,k): k\in \ss\}$,    $\triangle^c : = \ss\times \ss\setminus \triangle$ is the complement of $\triangle$,  and $H(x,y,k)$ is an appropriate $d \times d$
orthogonal matrix. Let us explain the coupling given in  \eqref{sgxkyl}.  
 When $k=l$ and $x\neq y$, we can  take
\begin{equation}\label{e:H(x,y,k)}
  H(x,y,k)= I-2\frac{(x-y)(x-y)^\top}{|x-y|^2};
\end{equation} 
this is the so-called coupling by reflection. When $k=l$ and $x=y$, the coupling in \eqref{sgxkyl} is   the so-called march coupling.
When $k\neq l$, we use 
the independent coupling in \eqref{sgxkyl}.
We refer to \cite{ChenLi-89} and \cite{LindR-86} for  details about  different couplings.

We now construct a coupling process $(X(t),\La(t),Y(t), \Psi(t))$ as follows. Let $(X(t),Y(t))$ satisfy the following SDE in $\rr^{2d}$,
\beq{XY}\d \begin{pmatrix}X(t)\\
Y(t) \end{pmatrix}= \tau \bigl(X(t),\La(t),Y(t), \Psi(t)\bigr)\d W(t)+
\begin{pmatrix} b \bigl(X(t),\La(t)\bigr)\\
b \bigl(Y(t), \Psi(t)\bigr)\end{pmatrix} \d t, \eeq where
$\tau(x,k,y,l)$ is the 
defined in \eqref{sgxkyl} and $W(t)$ is an
$\rr^{2d}$-valued standard Brownian motion. Meanwhile, let $(\La(t), \Psi(t))$
be a discrete random process with space $\ss \times
\ss$ such that \begin{align}\label{LaLaprime}
\nonumber &\bP\{(\La(t\!+\!\delta), \Psi (\!t+\!\delta))\!=\!(m,n)|
(\La(t), \Psi (t))\!=\!(k,l), (X_t,Y_t)\!=\!(\phi,\psi)\}\\
&\quad =\begin{cases}q_{(k,l)(m,n)}(\phi,\psi)\delta +o(\delta),
& \hbox{if} \ (m,n)\ne (k,l), \\
1+q_{(k,l)(k,l)}(\phi,\psi)\delta +o(\delta), & \hbox{if} \ (m,n)=(k,l)
\end{cases}
\end{align} provided $\delta \downarrow 0$, where 
  $\bigl(q_{(k,l)(m,n)}(\phi,\psi)\bigr)$ is the basic coupling of
  $\bigl(q_{kl}(\phi)\bigr)$ and $\bigl(q_{kl}(\psi)\bigr)$ derived from the basic coupling of \eqref{eq-Q(x)-coupling}.\comment{by the following relation:
\begin{equation}\label{Qxy}
\begin{split}Q(\phi,\psi) f(k,l)&:=\!\!
\sum_{m,n \in \ss}\!\! q_{(k,l)(m,n)}(\phi,\psi)\bigl(f(m,n)- f(k,l)\bigr)\\
&=\sum_{m\in \ss}
(q_{km} (\phi)-q_{lm} (\psi))^{+} \bigl(f(m,l)-f(k,l)\bigr)\\
&\  +\sum_{m\in \ss} (q_{lm} (\psi)-q_{km} (\phi))^{+} \bigl(f(k,m)-f(k,l)\bigr)\\
&\  +\sum_{m\in \ss} q_{km} (\phi)\wedge q_{lm}(\psi)\bigl(f(m,m)-f(k,l)\bigr)
\end{split}\eeq
for any function $f$ on $\ss \times \ss$; refer to \cite{Chen04}.}  
Note that for any $f: \ss\times \ss\mapsto \R$ and   $ (k, l) \notin \triangle$,  we can rewrite \eqref{eq-Q(x)-coupling}  as \begin{align*}Q& (\phi,\psi) f(k,l)\\&  := 
\sum_{m\in \ss\setminus\{k ,l\}}(q_{km} (\phi)-q_{lm} (\psi))^{+} \bigl(f(m,l)-f(k,l)\bigr) + (q_{kl} (\phi)-q_{ll} (\psi))^{+} \bigl(f(l,l)-f(k,l)\bigr)\\
&\quad  +\sum_{m\in \ss\setminus\{k ,l\}} (q_{lm} (\psi)-q_{km} (\phi))^{+} \bigl(f(k,m)-f(k,l)\bigr) +  (q_{lk} (\psi)-q_{kk} (\phi))^{+} \bigl(f(k,k)-f(k,l)\bigr) \\
& \quad +\sum_{m\in \ss\setminus\{k ,l\}} q_{km} (\phi)\wedge q_{lm}(\psi)\bigl(f(m,m)-f(k,l)\bigr)  \\ 
&\quad  +   q_{kk} (\phi)\wedge q_{lk}(\psi) \bigl(f(k,k)-f(k,l)\bigr)  +   q_{kl} (\phi)\wedge q_{ll}(\psi) \bigl(f(l,l)-f(k,l)\bigr) \\ 
 &  =  \sum_{m\in \ss\setminus\{k ,l\}} (q_{km} (\phi)-q_{lm} (\psi))^{+} \bigl(f(m,l)-f(k,l)\bigr)  \\ 
  & \quad  +\sum_{m\in \ss\setminus\{k ,l\}} (q_{lm} (\psi)-q_{km} (\phi))^{+} \bigl(f(k,m)-f(k,l)\bigr)      \\
 & \quad +  \sum_{m\in \ss\setminus\{k ,l\}} q_{km} (\phi)\wedge q_{lm}(\psi)\bigl(f(m,m)-f(k,l)\bigr) \\ 
 & \quad + q_{kl} (\phi) \bigl(f(l,l)-f(k,l)\bigr) + q_{lk} (\psi) \bigl(f(k,k)-f(k,l)\bigr).
\end{align*} This says that when  the coupling process $(\La,  \Psi)$ leaves  state $ (k, l)\notin \triangle$, it jumps to one of the  states in $\{(k, m), m\in \ss\setminus\{l\}\}\cup\{(m, l), m\in \ss\setminus\{k\}\}\cup\{(m, m), m\in \ss\} $ with respective   rates \begin{align}
\label{e-cp-jump-rates}
 \nonumber
q_{(k,l) (m, l)}(\phi,\psi) &= (q_{km} (\phi)-q_{lm} (\psi))^{+}, \ \quad  m \in \ss\setminus\{k ,l\}, \\
\nonumber q_{(k,l) (k, m)}(\phi,\psi) & = (q_{lm} (\psi)-q_{km} (\phi))^{+},  \ \quad m \in \ss\setminus\{k ,l\}, \\
q_{(k,l) (m, m)}(\phi,\psi) & =  q_{km} (\phi)\wedge q_{lm}(\psi),   \ \quad m \in \ss\setminus\{k ,l\}, \\
\nonumber q_{(k,l) (l, l)}(\phi,\psi) & = q_{kl} (\phi), \\
\nonumber q_{(k,l) (k, k)}(\phi,\psi) & =  q_{lk} (\psi).
\end{align}
 The sum of these rates is \begin{align}\label{e-cp-rate}
\nonumber -q_{(k,l)(k,l)}(\phi,\psi) & = \sum_{m \in \ss\setminus\{k ,l\}}  [ (q_{km} (\phi)-q_{lm} (\psi))^{+} +  (q_{lm} (\psi)-q_{km} (\phi))^{+} +  q_{km} (\phi)\wedge q_{lm}(\psi)] \\
 \nonumber & \quad + q_{kl} (\phi) +  q_{lk} (\psi) \\
 & = \sum_{m \in \ss\setminus\{k ,l\}}  q_{km} (\phi)\vee q_{lm}(\psi) +  q_{kl} (\phi) +  q_{lk} (\psi). 
\end{align} This is the rate at which the basic coupling process $(\La, \Psi)$ exits state $(k, l) \in \triangle^c$. Using the last three equations of \eqref{e-cp-jump-rates}, we see  that when $(\La, \Psi)$ exits  $(k, l) \in \triangle^c$, it enters $\triangle$ with rate \begin{align}\label{eq:diag_rate}
  \sum_{m\in \ss\setminus\{k, l\}}  q_{km} (\phi)\wedge q_{lm}(\psi) + q_{kl}(\phi) + q_{lk}(\psi). 
\end{align}

Using a similar argument as that for \eqref{e-cp-jump-rates}, when the coupling process $(\La,  \Psi)$ leaves  state $ (k, k)\in \triangle$, it jumps to  one of the states in $\{(k, m), (m, k), \text{ or }(m, m), m\in \ss\setminus\{k\}\} $ with respective   rates
\begin{align*} 
  q_{(k,k) (m, k)} (\phi,\psi)& = (q_{km} (\phi)-q_{km} (\psi))^{+},  \ \quad m \in \ss\setminus\{k \}, \\
 q_{(k,k) (k, m)}(\phi,\psi) & = (q_{km} (\psi)-q_{km} (\phi))^{+},  \ \quad m \in \ss\setminus\{k \}, \\
  q_{(k,k) (m, m)}(\phi,\psi) & =  q_{km} (\phi)\wedge q_{km}(\psi),   \ \quad m \in \ss\setminus\{k\}.
\end{align*} The sum of these rates is 
\begin{align}
\label{e2-cp-rate}
\nonumber   -q_{(k,k)(k,k)} (\phi,\psi)& = \sum_{m\in  \ss\setminus\{k \}} [(q_{km} (\phi)-q_{km} (\psi))^{+} +  (q_{km} (\psi)-q_{km} (\phi))^{+} + q_{km} (\phi)\wedge q_{km}(\psi)]  \\
    &  = \sum_{m\in  \ss\setminus\{k \}} q_{km} (\phi)\vee q_{km} (\psi).
\end{align}

Using the interlacing procedure together with exponential killing method as used before,
we can prove that for any given 
initial data $(\phi,k,\psi,l)\in \C\times \ss\times \C\times \ss$, there exists a coupling process $(X(t),\La(t),Y(t), \Psi(t))$ satisfying  the system (\ref{XY})--(\ref{LaLaprime}) with   $(X(t),Y(t))=(\phi(t),\psi(t))$ for $t\in [-r,0]$ and $(\Lambda(0), \Psi(0))=(k,l)$.
Moreover, $(X_t,\Lambda(t),Y_t, \Psi(t))$ is a Markov-Feller process. 

We can define the operator associated with  the coupling process $(X(t),\La(t),Y(t), \Psi(t))$ of \eqref{XY}--\eqref{LaLaprime} as follows. First, we define the matrix 
\begin{align}\label{e:A(xkyl)}
 \nonumber \wdh a(x,k, y, l) & := \tau(x,k,y,l)\tau(x,k,y,l)^\top  \\
 \nonumber & =  {\mathbf{1}}_\triangle (k,l){\mathbf{1}}_{\{x\ne y\}}\begin{pmatrix}
    a(x,k) & \sigma(x,k)^\top H(x,y,k)\sigma(y, k)\\ \sigma(y, k) H(x,y,k) \sigma(x,k)^\top & a(y,k)
  \end{pmatrix}\\ 
 \nonumber    & \quad +{\mathbf{1}}_\triangle (k,l){\mathbf{1}}_{\{x=y\}}\begin{pmatrix}
    a(x,k) & \sigma(x,k)^\top \sigma(y,k)\\ \sigma(y,k) \sigma(x,k)^\top & a(y,k)
  \end{pmatrix}\\ 
  & \quad +{\mathbf{1}}_{\triangle^c}(k,l)\begin{pmatrix}
    a(x,k) & 0\\
    0 & a(y,l)
  \end{pmatrix},
\end{align} 
where   $a(x,k)=\sg(x,k)\sg(x,k)^\top$. Similar to the definition of \eqref{e:op-basic-coupling}, for any $f \in C_c^2(\R^d\times \ss\times \R^d\times \ss)$,  we define
\begin{align}\label{e:wdh-op-reflection}
  \wdh \op f(\phi,k,\psi,l) & =    \wdh \varOmega_{\text{d}}   f(\phi(0),k, \psi(0),l) +   \wdt \varOmega_{\text{s}}f(\phi,k, \psi,l), \quad \forall (\phi, k, \psi, l) \in \C\times \ss\times \C\times \ss,
\end{align} where  $\wdt \varOmega_{\text{s}}f(\phi,i, \psi,j)$ is defined in \eqref{eq-Q(x)-coupling}, and 
\begin{align}\label{e:wdh-op}
  \wdh \varOmega_{\text{d}} f(\phi,k,\psi,l) & = \frac12 \tr(\wdh a(x, k, y, l) D^2 f(x,k,y, l)) + \lan b(x, k, y, l), Df(x,k,y,l)\ran,
\end{align} in which  $b(x, k, y, l)$ is defined in \eqref{eq:double-b}.

Let $\bP^{(\phi,k,\psi,l)}$ denote the distribution of the coupling process 
$(X_t,\La(t),Y_t, \Psi(t))$ starting from $(\phi,k,\psi,l)$, and $\E^{(\phi,k,\psi,l)}$ denote the expectation with respect to  $\bP^{(\phi,k,\psi,l)}$.
Furthermore, we set 
 the {\it coupling time} of $(X_t,\La(t))$ and $(Y_t, \Psi(t))$ as
\begin{equation}\label{eq:T}
  T := \inf \{t \ge 0: X_t=Y_t,\ \La(t)= \Psi(t)\}.
\end{equation}  We also  define   the {\it meeting time} of $(X(t),\La(t))$ and $(Y(t), \Psi(t))$ as
\begin{equation}\label{eq:hatT}
  \widehat{T} = \inf \{t \ge 0: X(t)=Y(t),\ \La(t)= \Psi(t)\}.
\end{equation} Obviously, the meeting time does not exceed the coupling time, i.e., it is always true that $\widehat{T}\le T$.  As alluded in the introduction,  the first meeting time $\wdh T  < \infty$ does not guarantee a successful coupling for the processes $(X_t,\La(t))$ and $(Y_t, \Psi(t))$. Because of the past-dependent transition rates for $(\La,  \Psi)$, the processes may still separate from each other after $\wdh T$. This is in sharp contrast to the classical regime-switching diffusions without delay.

\begin{Definition}\label{defsuccess} The   coupling $(X_t,\La(t), Y_t, \Psi(t))$ is said to be {\em successful} if
$$\bP^{(\phi,k,\psi,l)}\bigl(T<\infty\bigr)=1, \ \ \forall (\phi,k)\ne (\psi,l).$$\end{Definition}


In order to 
show that the coupling $(X_t,\La(t), Y_t, \Psi(t))$ defined in \eqref{XY}--\eqref{LaLaprime} is successful, for each $k\in \ss$, we first  consider the coupling of $X^{(k)}(t)$ and itself.
To this end, in view of \eqref{sgxkyl}, for any $(x,y,k) \in \rr^d \times \rr^d \times \ss$, set a $2d \times 2d$ matrix as follows:
$$\tau(x,y,k): =\tau(x, k, y, k) = {\mathbf{1}}_{\{x\ne y\}}
\begin{pmatrix}
  \sg(x,k) & 0\\
  \sg(y,k)H(x,y,k) & 0
\end{pmatrix} 
+{\mathbf{1}}_{\{x=y\}}
 \begin{pmatrix}
  \sg(x,k) & 0\\
  \sg(x,k) & 0
\end{pmatrix} ,$$
where $H(x,y,k)$ is an appropriate $d \times d$ orthogonal matrix. For example, we  can take  $H(x,y,k)$ as in \eqref{e:H(x,y,k)}, yielding the reflection coupling.

Let $(X^{(k)}(t),Y^{(k)}(t))$ satisfy the following SDE in $\rr^{2d}$
\beq{XkYk}\d \begin{pmatrix} X^{(k)}(t)\\
Y^{(k)}(t)
 \end{pmatrix}  = \tau \bigl(X^{(k)}(t),Y^{(k)}(t),k\bigr)\d W(t)+
 \begin{pmatrix} b \bigl(X^{(k)}(t),k\bigr)\\
b \bigl(Y^{(k)}(t),k\bigr) \end{pmatrix} \d t. \eeq
Then we set the coupling time of $X^{(k)}(t)$ and $Y^{(k)}(t)$ as follows:
\begin{equation}\label{eq:Tk}
  T^{(k)} = \inf \{t \ge 0: X^{(k)}(t)=Y^{(k)}(t)\}.
\end{equation} 

\begin{Assumption}\label{Tk}{\rm \begin{itemize}
  \item[(i)] There exists a positive constant  $M$ such that  
the coupling time $T^{(k)}$ defined in \eqref{eq:Tk} satisfies
\begin{equation}\label{supTk}\bP^{(x,y,k)} \bigl(T^{(k)}<M \bigr)\ge \frac{1}{2} \ \
\hbox{for all} \ (x,y) \in \rr^{2d} \ \hbox{and} \ k\in \ss,\end{equation} where $\bP^{(x,y,k)}$ denotes the distribution of the coupling 
$(X^{(k)}(t),Y^{(k)}(t))$ of \eqref{XkYk} with initial condition   $(x,y)$.
 
\item[(ii)] There exists a positive constant $\alpha>0$ such that for any $(\phi,\psi)\in \C \times \C$
and $(k,l)\in \ss^{2}\setminus \triangle$, \begin{equation}\label{qh}\sum_{m\in
\ss\setminus \{k,l\}}q_{km}(\phi)\wedge q_{lm}(\psi)+q_{kl}(\phi)+q_{lk}(\psi)\ge
\alpha>0. \end{equation}\end{itemize}
}\end{Assumption}

\begin{Remark}\label{onTk}  1) Assumption~\ref{Tk} (i) will be satisfied if for every $k\in \ss$, we can construct a successful coupling
$(X^{(k)}(t),Y^{(k)}(t))$ of $X^{(k)}(t)$ and $Y^{(k)}(t)$, and moreover the moment of the corresponding coupling time $T^{(k)}$ is
uniformly bounded with respect to the starting points. 

Let us present a sufficient condition for the uniform boundedness of $\E^{(x, y,k)}[T^{(k)}]$, where   $\E^{(x, y,k)}$ denotes the expectation with respect to the probability $\P^{(x,y,k)}$.

 Suppose there exists  a nonnegative and uniformly bounded function $F\in C^2([0, \infty))$ satisfying   \begin{align}\label{F-condition}
  \mathcal L^{(k)} F( |x- y| )\le -K < 0, \quad \text{for all } (x, y) \in \R^d\times \R^d \text{ with } x\neq y,
 \end{align}   where   $\mathcal L^{(k)}$ is the operator associated with the process \eqref{XkYk} and  $K$ is a positive constant,  then we have \begin{equation}\label{eq:Tk-moment}
  \sup_{(x,y)\in \rr^{2d}} \E^{(x,y,k)}[T^{(k)}] < \infty.
 \end{equation} Indeed, \eqref{eq:Tk-moment} Obviously holds for $x =y$. For $x\neq y$,   we can use the It\^o's formula and \eqref{F-condition} to derive that for any $t\ge 0$,
$$\E^{(x,y,k)}[F(|X^{(k)}(t\wedge T^{(k)})-Y^{(k)}(t\wedge T^{(k)})|)] \le F(|x-y|) - K \E^{(x,y,k)}[t\wedge T^{(k)}].$$
   This, in turn, implies that $$\E^{(x,y,k)}[T^{(k)}] \le \frac{F(|x-y|)}{K} \le \frac{\|F\|_\infty}{K} < \infty. $$ Thus \eqref{eq:Tk-moment} holds. Alternatively, one can use impose conditions on the coefficients of \eqref{eq:X} so that \eqref{eq:Tk-moment} holds;  we refer to \cite{Mao06,ShaoX-13} for details.
 
2) Assumption~\ref{Tk} (ii)  is a technical condition, and it guarantees that for the coupling constructed by
\eqref{eq-Q(x)-coupling}, $(\La(t), \Psi(t))$ can transition  from $  \triangle^c$ to $\triangle$ in finite time a.s.; see Proposition \ref{varsigmafinite}.  In view of   \eqref{eq:diag_rate}, Assumption~\ref{Tk} (ii) implies that the rate at which $(\La, \Psi)$ enters the diagonal $\triangle$ upon   exiting any state $(k, l) \in \triangle^c$ is uniformly bounded from below. Intuitively,  this condition governs the rate at which the two components of $(\La, \Psi)$ are able to coalesce, while \eqref{(1.6)} in Assumption~\ref{I1} controls the duration they may stay together after coalescence. Together with Assumption~\ref{Tk}~(i), these conditions ensure the successful coupling of the segment process $(X_t, \La(t))$.
\end{Remark}

\begin{Theorem}\label{success1}
Suppose Assumptions~\ref{I1} and \ref{Tk} hold. Then the coupling $(X_t,\La(t),Y_t, \Psi(t))$
constructed by   \eqref{XY} and  \eqref{LaLaprime} is successful, that is
\beq{Tfinite}\bP^{(\phi,k,\psi,l)} (T <\infty) =1, \qquad \text{for any }\ (\phi,k) \ne (\psi,l).\eeq
\end{Theorem}
Before we present the proof of   Theorem \ref{success1} in the end of the section, let us first state and prove the following corollary, which indicates that the process $(X_t, \La(t))$ is ergodic with an invariant measure $\pi$. 

\begin{Corollary}\label{success3}
Suppose Assumptions~\ref{I1}, \ref{assump-Q}, and \ref{Tk} hold. Then 
 the Markov-Feller process $(X_t,\La(t))$ has an invariant measure $\pi\in\mathcal{P}(E)$, and its
transition probability $P(t,(\phi,k), \cdot)$ converges to $\pi$ in total variation norm  as $t \to \infty $ for every $(\phi,k)\in \C \times \ss$.
\end{Corollary}

\para{Proof.} 
By virtue of the coupling inequality, for any $u, t>0$ we obtain
\begin{eqnarray*}
	 &&\|P(t+u, (\phi,k), \cdot)
	-P(t, (\phi,k), \cdot) \|_{\TV} \\
	&&\quad =\bigg\|P(t,(\phi,k), \cdot) - \displaystyle\int_{E}  P(u,(\phi,k), \d \psi \times\d\{l\})P(t,(\psi,l), \cdot)\bigg\|_{\TV} \\
	&&\quad \le  \displaystyle\int_{E} P(u,(\phi,k), \d \psi \times\d\{l\})
	\|P(t, (\phi,k), \cdot)-P(t, (\psi,l), \cdot) \|_{\TV} \\
	&&\quad\le 2
	\displaystyle\int_{E} P(u,(\phi,k), \d \psi \times\d\{l\}) \bP^{(\phi,k,\psi,l)} (T> t)
	\to 0, \quad \hbox{as} \quad t \to \infty. 
\end{eqnarray*}
This shows that the family $\{P(t, (\phi, k), \cdot)\}_{t\geq0}$ is Cauchy in the total variation norm and hence tight in $(\mathcal P(E), \|\cdot\|_{\TV})$ thanks to Lemma \ref{lem-tight} below. Combining the Feller property of the process $(X_t, \Lambda(t))$ with the Krylov-Bogoliubov theorem \cite[Corollary 3.1.2]{DaPratoZ-96}, we establish the existence of an invariant probability measure $\pi \in \mathcal{P}(E)$.
	
Similarly, applying the coupling inequality again yields
\begin{align*}  \|P(t,(\phi,k), \cdot)
-\pi (\cdot) \|_{\TV}  
&  =\bigg\|P(t,(\phi,k), \cdot) - \displaystyle\int_{E} \d \pi
(\psi,l) P(t,(\psi,l), \cdot) \bigg\|_{\TV} \\
&  \le  \displaystyle\int_{E} \d \pi (\psi,l)
\|P(t, (\phi,k), \cdot)-P(t, (\psi,l), \cdot) \|_{\TV} \\
& \le 2  \displaystyle\int_{E} \d \pi (\psi,l)
 \P^{(\phi,k,\psi,l)} (T> t)
\to 0, \quad \hbox{as} \quad t \to \infty. \end{align*} This shows that   $(X_t,\La(t))$ is ergodic and completes the proof. \qed

\begin{Lemma}\label{lem-tight}
		Let $\{\mu_n\}$ be a sequence of probability measures on $\mathcal{B}(E)$ that is Cauchy in the total variation norm. Then $\{\mu_n\}$ is tight. 
\end{Lemma} 
\para{Proof.} Since $\{\mu_n\}$ is Cauchy in the total variation norm, for any $\e > 0$, there exists an $N \in \mathbb{N}$ so that \begin{align*}
  \|\mu_n- \mu_N\|_\TV  < \frac\e2, \quad \forall n \ge N.
\end{align*} On the other hand, since $E = \C\times \ss$ is a Polish space, for each $i =1, \dots, N$, we can find compact subsets $K_i\subset E$ so that 
\begin{align*}
  \mu_i(K_i^c) < \frac\e2, \quad \forall i =1, \dots, N.
\end{align*} Now let $K : = \bigcup_{i=1}^N K_i$, which is still a compact subset of $E$. We have \begin{align*}
  \mu_i(K^c) \le \mu_i(K_i^c) < \frac\e 2, \quad \forall i =1, 2,\dots, N,
\end{align*} and for all $n \ge N$
\begin{align*}
  \mu_n(K^c) \le \mu_n(K_N^c) \le \mu_N(K_N^c) + \|\mu_n- \mu_N\|_\TV < \e.
\end{align*} Therefore $\{\mu_n\}$ is tight. \qed
		
		

Clearly, the  proof of Corollary \ref{success3} relies on the fact that  the coupling $(X_t,\La(t),Y_t, \Psi(t))$ is successful. To proceed, let us introduce some notation.
Let $\{\eta_{m}\}_{m\ge 0}$ be a sequence of stopping times defined
by $\eta_0 =0,$ and for $m \ge 1$,   \beq{eta0m}   
\eta_m =\inf \{t\ge \eta_{m-1}:
(\La(t), \Psi(t))\ne
(\La(\eta_{m-1}), \Psi(\eta_{m-1}))\}.\eeq Let
$\{\zeta_{m}\}_{m\ge 1}$ be a sequence of stopping times defined by
\beq{devarsigma1}\zeta_1 =\inf \{t\ge 0: (\La(t), \Psi(t))\in
\triangle\},\eeq and for $n=1,2,\cdots$,
\beq{devarsigma2}\barray\zeta_{2n} \ad=\inf \{t\ge \zeta_{2n-1}:
(\La(t), \Psi(t))\in \ss^{2}\setminus\triangle\},\\
\zeta_{2n+1} \ad=\inf \{t\ge \zeta_{2n}: (\La(t), \Psi(t))\in
\triangle\}.\earray\eeq That is, $\{\zeta_{2n-1}\}_{n\ge 1}$ is the
sequence of successive hitting times to $\triangle$, whereas
$\{\zeta_{2n}\}_{n\ge 1}$ is the sequence of successive exit times from
$\triangle$. Clearly, the sequence $\{\zeta_{m}\}_{m\ge 1}$ is a
subsequence of $\{\eta_{m}\}_{m\ge 0}$. 
The sequence $\{\zeta_{m}\}_{m\ge 1}$ leads to 
 a sequence of
temporal intervals. On $[\zeta_{2n-1},\zeta_{2n})$, $\La(t)$ and $ \Psi(t)$
coincide; while on $[\zeta_{2n}, \zeta_{2n+1})$, $\La(t)$ and
$ \Psi(t)$ are different. 
For the coupling
$(X_t,\La(t),Y_t, \Psi(t))$,  
Proposition~\ref{varsigmafinite} below shows that
the first hitting time $\zeta_{1} $ is finite almost surely with respect
to $\bP^{(\phi,k,\psi,l)}$ for any  $(\phi,\psi)\in \C \times \C$ and $(k,l)\in
\ss^{2}\setminus \triangle$. This  fact will play an important roll in the proof of Theorem \ref{success1}.


\begin{Proposition}\label{varsigmafinite}
Suppose Assumptions~\ref{I1} and ~\ref{Tk} (ii) hold, then  we have \beq{eqvarsigmafinite}\bP^{(\phi,k,\psi,l)}
(\zeta_{1} <\infty) =1, \quad  \ \text { for all }\  (\phi,\psi)\in \C \times \C  \ \text{ and   }
\ (k,l)\in \ss^{2}\setminus \triangle. \eeq
\end{Proposition}

\para{Proof.} We consider the function $f(x, k, y, l):=\1_{ \{k\neq l\}} $,   $(x, k, y, l)\in \R^d\times \ss \times \R^d\times \ss$. Straightforward calculations reveal that for the operator $\wdh\op$ defined in \eqref{e:wdh-op-reflection}, we have \begin{align*}
  \wdh \op f(\phi, k, \psi, l) & =  \wdt \varOmega_{\text{s}}f(\phi, k, \psi, l)\\
   & = - \sum_{m\in \ss\setminus\{k, l\}} q_{km} (\phi) \wedge q_{lm}(\psi)- q_{kl} (\phi) - q_{lk}(\psi) 
   \\ & \le -\alpha,  
\end{align*} for all $(\phi, k, \psi, l)\in \C\times \ss\times \C\times \ss$ with $k\neq l$, where the last inequality follows from   Assumption~\ref{Tk} (ii). Then it follows that  
 for any $(\phi,k,\psi, l)\in \C \times \ss\times \C\times \ss$ with $k\neq l$, we have  \begin{align*}
\E^{{(\phi,k,\psi,l)}}& [  f(X(t\wedge\zeta_1), \La(t\wedge\zeta_1), Y(t\wedge\zeta_1),  \Psi(t\wedge\zeta_1))]\\ & = f(X(0), \La(0), Y(0),  \Psi(0)) + \E^{{(\phi,k,\psi,l)}}\bigg[  \int_0^{t\wedge\zeta_1} \wdh \op f(X(s), \La(s), Y(s),  \Psi(s)) ds\bigg]\\ & \le f(X(0), \La(0), Y(0),  \Psi(0)) - \alpha \E^{{(\phi,k,\psi,l)}}[t\wedge\zeta_1].
\end{align*} Using the definition of $f$, we can rewrite the above equation as $\alpha \E^{{(\phi,k,\psi,l)}}[t\wedge\zeta_1] \le 1$. Letting $t\to\infty$ yields $\E^{{(\phi,k,\psi,l)}}[\zeta_1] \le \frac1\alpha$. Note that this estimate holds for all  $(\phi,\psi)\in \C\times \C$ and   all $(k, l)\in \ss^2$. In other words, we have the uniform bound  \begin{equation}\label{e:zeta1mean}
 \sup_{(\phi,k,\psi,l) \in \C\times \ss\times \C\times \ss} \E^{{(\phi,k,\psi,l)}}[\zeta_1] \le \frac1\alpha.
\end{equation} This bound certainly implies  \eqref{eqvarsigmafinite}.
\qed

With Proposition~\ref{varsigmafinite} established, we  are now ready to  prove Theorem~\ref{success1}. 
\para{Proof of Theorem \ref{success1}.} \ 
To prove \eqref{Tfinite}, we need to give some estimation on the
first exit time $\zeta_{2}$ via the first jump time $\eta_1$. In
view of the definitions $\eta_1$, $\zeta_{1}$ and $\zeta_{2}$
in (\ref{eta0m}), (\ref{devarsigma1}), and (\ref{devarsigma2}), respectively,
when $(X_0,\La(0),Y_0, \Psi(0))=(\phi,k,\psi,k)$, we have that $\zeta_{1}=0$ and $\zeta_{2}\ge
\eta_1$. Then it follows   that
$$\{T \in [0, \zeta_2)\} \supset \{T\in [0,\eta_1)\}\supset \{\widehat{T}\in [0,\eta_1), \widehat{T}+r\le \eta_1\}
\supset \{\widehat{T}< M, \eta_1 \ge M+r\},$$ where $M$ is the constant given in Assumption~\ref{Tk} (i).
Thus
\begin{align*}
  \bP^{(\phi,k,\psi,k)}\bigl(T\in[0,\zeta_2)\bigr) & \ge \bP^{(\phi,k,\psi,k)}\bigl(T\in [0,\eta_1)\bigr)\\ &  \ge
\bP^{(\phi,k,\psi,k)}\bigl(\eta_1 \ge 
 M+r\bigr)\bP^{(\phi,k,\psi,k)}\bigl(\widehat{T}<M|\eta_1 \ge  M+r\bigr).
\end{align*}
To estimate 
 the probabilities in the above equation,  we consider the element $q_{(k,k)(k,k)}(\phi,\psi)$ in the
coupling matrix determined by \eqref{eq-Q(x)-coupling}. We know from \eqref{(1.6)} and 
(\ref{e2-cp-rate}) that for any $(\phi,\psi)\in \C \times \C$ and $k\in \ss$,
$$-q_{(k,k)(k,k)}(\phi,\psi)=\sum_{m\in \ss\setminus\{k\}}q_{km}(\phi)\vee
q_{km}(\psi)\le H.$$ Hence, using the proofs of \cite[Lemmas 3.2
and 3.3]{XiZ-06}  or applying the  last displayed equation on p. 103 of \cite{Skorohod-89}, we derive that 
\begin{equation}
\label{e-eta1-distribution}
\bP^{(\phi,k,\psi,k)}\bigl(\eta_1 \ge t
\bigr) \ge \exp \bigl\{-H t\bigr\}, \quad \text{ for all }t >0.
\end{equation} 
On the other hand,  conditional on  $\{\eta_1 \ge M+r\}$, the evolution of the process $(X(t),Y(t))$ is the 
same as the that  of the coupling process $(X^{(k)}(t),Y^{(k)}(t))$ on 
the temporal interval $[0,M+r]$.  Therefore, we have \begin{equation}
\label{e-T-hat}
\bP^{(\phi,k,\psi,k)}\bigl(\widehat{T}< M|\eta_1 \ge M+ r\bigr)=
\bP^{(\phi(0),\psi(0), k)}\bigl(T^{(k)}< M\bigr).
\end{equation}
Thus, it follows from \eqref{e-eta1-distribution}, \eqref{e-T-hat}, and   Assumption~\ref{Tk} (i)  
that for any   $(\phi, \psi) \in \C\times \C$
and $k\in \ss$,
\begin{align}\label{T2} \nonumber   \bP^{(\phi,k,\psi,k)}\bigl(T\in [0,\zeta_2)\bigr) &  \ge   \bP^{(\phi,k,\psi,k)}\bigl(T\in [0,\eta_1)\bigr)\\ 
\nonumber  &\ge \bP^{(\phi,k,\psi,k)}\bigl(\eta_1 \ge M+r\bigr)\bP^{(\phi,k,\psi,k)}\bigl(\widehat{T}< M|\eta_1 \ge M+r\bigr)\\
\nonumber  &\ge \bP^{(\phi,k,\psi,k)}\bigl(\eta_1 \ge
M+r\bigr)\bP^{(\phi(0),\psi(0),k)}\bigl(T^{(k)}< M\bigr)\\
& \ge \frac{1}{2}\exp\{-H(M+r)\}=:\delta_2>0,
\end{align} 
note that  the positive constant $\delta_2$ is independent of the initial data $\phi$,
$\psi$ and $k$.  In other words, \eqref{T2} holds uniformly for all $\phi,\psi \in \C$ and $k \in \ss$. 

 Recall from Proposition~\ref{varsigmafinite} that 
$\zeta_{1} <\infty$ a.s. with respect to $\bP^{(\phi,k,\psi,l)}$, where $k\neq l$. Note also that $\La(\zeta_1) =  \Psi(\zeta_1)$. 
Therefore we can apply  the strong Markov
property and \eqref{T2} to derive \beq{T1}\begin{aligned}
\bP^{(\phi,k,\psi,l)}\bigl(T\notin
[\zeta_{1},\zeta_{2})\bigr) & =\E^{(\phi,k,\psi,l)}
\big[\bP^{(X_{\zeta_{1}},\La(\zeta_{1}),Y_{\zeta_{1}}, \Psi(\zeta_{1}))}
\bigl(T\notin [0,\zeta_{2})\bigr)\big]\\
& \le\E^{(\phi,k,\psi,l)}[1-\delta_2] = 1-\delta_2.
 \end{aligned} \eeq
 Then it follows from the
strong Markov property that for any positive integer $m$,
\begin{align*}
&\bP^{(\phi,k,\psi,l)}\Bigg\{T\notin
\bigcup_{n= 1}^{m}[\zeta_{2n-1},\zeta_{2n})\Bigg\}\\
&= \E^{(\phi,k,\psi,l)}\Bigg[\bigl\{T\notin
\bigcup_{n= 1}^{m-1}[\zeta_{2n-1},\zeta_{2n})\bigr\};
\bP^{(X_{\zeta_{2m\!-\!2}},\La(\zeta_{2m\!-\!2}),Y_{\zeta_{2m\!-\!2}}, \Psi(\zeta_{2m\!-\!2}))}
\bigl(T\notin   [\zeta_{2m-1},\zeta_{2m})\bigr) \Bigg] \\
&\leq  \bigl(1-\delta_2\bigr) \bP^{(\phi,k,\psi,l)}\Bigg\{T\notin
\bigcup_{n= 1}^{m-1}[\zeta_{2n-1},\zeta_{2n})\Bigg\}\\ 
&\le \bigl(1-\delta_2\bigr)^{m}.
\end{align*} 
Therefore we have 
\begin{align*}
\bP^{(\phi,k,\psi,l)}\bigl(T=\infty\bigr) &   \le
\bP^{(\phi,k,\psi,l)}\Bigg\{T\notin \bigcup_{n=
1}^{\infty}[\zeta_{2n-1},\zeta_{2n})\Bigg\}\\ & = \lim_{m\to
\infty}\bP^{(\phi,k,\psi,l)}\Bigg\{T\notin \bigcup_{n=
1}^{m}[\zeta_{2n-1},\zeta_{2n})\Bigg\}\le \lim_{m\to \infty}
\bigl(1-\delta_2\bigr)^{m}=0;
\end{align*}
this clearly implies
(\ref{Tfinite}). The proof of Theorem \ref{success1} is complete. \qed

\begin{Remark}
	It is worth noting that the proof of Theorem \ref{success1} does not depend on the specific form of the coupling in \eqref{XY}. The only key requirement is that once the marginal processes coalesce, they move together afterwards.  
\end{Remark}

\section{Strong Ergodicity}\label{sect-stability}
 We have shown in Corollary \ref{success3} that under Assumptions~\ref{I1}, \ref{assump-Q}, and \ref{Tk}, the process $(X_t,\La(t))$ is ergodic with an invariant measure $\pi$ and $\| P(t, (\phi, k), \cdot)- \pi(\cdot)\|_\TV \to 0 $ as $t\to\infty$.  Theorem \ref{thm:stroergo} below further demonstrates that under these conditions, the process $(X_t,\La(t))$ is strongly ergodic with an explicit exponential convergence rate that is uniformly bounded from below.

For Markov processes, the strong ergodicity (or uniform ergodicity) in
the sense of convergence in variation norm is the strongest one; see
\cite{Chen04, MeynT-09}.  Strong ergodicity for certain Markov processes has been studied using coupling methods. This includes work on one-dimensional diffusion processes and birth-death processes \cite{Mao02, Mao06}. See also \cite{XiYin-09} which investigates  strong ergodicity   for a mean-field model with a continuous-state-dependent switching process without delay.  

To proceed, we first recall the  
definition of strong ergodicity.

\begin{Definition}{\cite[Definition 1.1]{Mao06}}
	The  Markov process $(X_t,\La(t))$ is said to be {\em strongly ergodic in the variational norm} (in short, strongly ergodic) if there exists an $\varepsilon>0$ such that
	\beq{defstre}
		\sup\{\|\bP(t, (\phi, k), \cdot)-\pi (\cdot)\|_{\TV}: (\phi, k)\in \C\times \ss\}=O({ \rm e}^{-\varepsilon t}) \quad\text{  as }\ t\rightarrow\infty,
	\eeq where $\pi$ is the invariant distribution of the process $(X_t,\La(t))$.    

	For $\gamma\geq2$, define
	\beq{hatalpha}
	\beta(\gamma)=\sup\{\varepsilon\geq0:\sup\{\|\bP(t, (\phi, k), \cdot)-\pi (\cdot)\|_{\TV}: (\phi, k)\in \C\times \ss\}\leq\gamma{ \rm e}^{-\varepsilon t} \text{ for all $t\geq0$}\}.
	\eeq  Obviously, $\beta(\gamma)$ is increasing with respect to $\gamma$. Therefore we can define 
	  $\beta:=\beta(\infty)=\lim_{\gamma\rightarrow\infty}\beta(\gamma)$.
\end{Definition}

\begin{Theorem}\label{thm:stroergo} 
	Suppose Assumptions~\ref{I1}, \ref{assump-Q},  
	and \ref{Tk} hold. Then the Markov process $(X_t,\La(t))$ is strongly  ergodic. In particular, for any $\gamma \ge  2$, defining \begin{align}\label{eq:lambda_defn}
     \lambda(\gamma) :=  \frac{1-\frac2\gamma}{4 e^{H(M+r)} (M+r + \frac{2}{ \alpha} + 1)}, \end{align} where $r >0$ is the delay duration in \eqref{eq:Lambda}, and the positive constants $H, M$, and $\alpha $ are respectively given in \eqref{(1.6)},  \eqref{supTk}, and \eqref{qh}, 
     we have 
\begin{align}\label{eq:uni-erg}
  \sup_{(\phi, k)\in \C\times \ss} \|P(t, (\phi, k), \cdot)- \pi(\cdot)\|_\TV \le \gamma e^{-\lambda(\gamma) t}, \quad \forall t\ge 0. 
\end{align} Therefore, 
	\begin{equation}\label{stre}  \beta(\gamma) \ge  \lambda(\gamma),  \ \forall \gamma > 2,  \quad  \text{ and } \quad 
	\beta\geq 
  \frac{1}{4 e^{H(M+r)} (M+r + \frac{2}{ \alpha} + 1)}>0.
  \end{equation}
\end{Theorem}
\para{Proof.}  We only need to prove \eqref{eq:uni-erg} because \eqref{stre} then follows from the definitions of $\beta(\gamma)$ and $\beta$.
Adopting the notation from the preceding discussion, for any $(\phi, \psi)\in\mathcal C\times\mathcal C$ and $k\in\mathbb S$, we deduce from part (i) of Assumption~\ref{Tk}  that
\begin{align}\label{kk} 
  \nonumber \bP^{(\phi,k,\psi,k)}(T<M+r) &\ge \bP^{(\phi,k,\psi,k)}(T<M+r, \eta_{1}\geq M+r)\\
	 \nonumber  &= \bP^{(\phi,k,\psi,k)}(\eta_{1}\geq M+r)\bP^{(\phi,k,\psi,k)}(T<M+r| \eta_{1}\geq M+r)\\
	 \nonumber  & \ge  \bP^{(\phi,k,\psi,k)}(\eta_{1}\geq M+r)\bP^{(\phi,k,\psi,k)}(\widehat{T}<M| \eta_{1}\geq M+r)\\
	 \nonumber  &=  \bP^{(\phi,k,\psi,k)}(\eta_{1}\geq M+r)\bP^{(\phi(0),\psi(0),k)}\bigl(T^{(k)}< M\bigr)\\
	& \ge \frac{1}{2}\exp\{-H(M+r)\}=\delta_2>0; 
\end{align} where the last equality follows from the definition of $\delta_2$ given in  \eqref{T2}.
On the other hand, with $N = N(\alpha): = \frac2\alpha +1 > 0$,  it follows from \eqref{e:zeta1mean} and the Markov inequality that 
\begin{align*}
  \P^{{(\phi,k,\psi,l)}}[\zeta_1\geq N] \le \frac{\sup_{(\phi,k,\psi,l) \in \C\times\ss\times \C\times \ss} \E^{(\phi,k,\psi,l)} [\zeta_1]}{N} \le \frac{\frac1\alpha}{\frac2\alpha+1} < \frac12. 
\end{align*} Therefore we have 
\beq{N}
\P^{{(\phi,k,\psi,l)}}[\zeta_1<N]\geq\frac{1}{2}.
\eeq
 For any initial conditions $(\phi, \psi)\in\mathcal C\times\mathcal C$ and distinct states $k\neq l\in \mathbb{S}$, an application of the strong Markov property combined with \eqref{kk} and \eqref{N} yields
\begin{align}\label{kneql}
  \nonumber  \bP^{(\phi,k,\psi,l)}(T<M+r+N) &\ge \bP^{(\phi,k,\psi,l)}(T<M+r+N, \zeta_{1}<N)\\
	\nonumber   &= \mathbb E^{(\phi,k,\psi,l)}\bigl[\mathbb E^{(\phi,k,\psi,l)}[\mathbf{1}_{\{T<M+r+N\}}\mathbf{1}_{\{\zeta_{1}<N\}}|\mathcal F_{\zeta_{1}}]\bigr]\\
	\nonumber   &= \mathbb E^{(\phi,k,\psi,l)}\bigl[\mathbf{1}_{\{\zeta_{1}<N\}}\mathbb E^{(\phi,k,\psi,l)}[\mathbf{1}_{\{T<M+r+N\}}|\mathcal F_{\zeta_{1}}]\bigr]\\
	\nonumber   &= \mathbb E^{(\phi,k,\psi,l)}\bigl[\mathbf{1}_{\{\zeta_{1}<N\}}\mathbb E^{(X_{\zeta_{1}},\Lambda(\zeta_{1}),Y_{\zeta_{1}}, \Psi(\zeta_{1}))}[\mathbf{1}_{\{T<M+r+N-\zeta_{1}\}}]\bigr]\\
	\nonumber   &\ge \mathbb E^{(\phi,k,\psi,l)}\bigl[\mathbf{1}_{\{\zeta_{1}<N\}}\mathbb E^{(X_{\zeta_{1}},\Lambda(\zeta_{1}),Y_{\zeta_{1}}, \Psi(\zeta_{1}))}[\mathbf{1}_{\{T<M+r\}}]\bigr]\\
	& \ge \delta_2\P^{{(\phi,k,\psi,l)}}[\zeta_1<N]  \ge \frac{\delta_{2}}{2}.
\end{align} 
On the other hand,  for any initial condition $(\phi, \psi)\in\mathcal C\times\mathcal C$ and $k\in\mathbb S$, we have from \eqref{kk} that 
\beq{keql}
\bP^{(\phi,k,\psi,k)}(T<M+r+N)\geq\bP^{(\phi,k,\psi,k)}(T<M+r)\geq\delta_{2}>\frac{\delta_{2}}{2}.
\eeq
By combining   \eqref{kneql} and \eqref{keql}, we conclude that for any $(\phi, \psi)\in\mathcal C\times\mathcal C$ and $(k, l)\in\mathbb S\times\mathbb S$, the following uniform estimate holds:
\beq{klleq}
 \bP^{(\phi,k,\psi,l)}(T<M+r+N)\geq\frac{\delta_{2}}{2}>0.   
\eeq In other words, for any $(\phi, k, \psi, l)\in\mathcal C\times\ss\times \mathcal C\times \ss$, we have
\beq{kl}
\bP^{(\phi,k,\psi,l)}\bigl(T  \geq R \bigr)\leq \rho,  
\eeq where $R: =  {M+r+N}$ and $\rho: = 1-\frac{\delta_2}{2} < 1$.
Using the Markov property and an inductive argument, we establish that for any $n \in \mathbb N$, 
\begin{align}\label{induction} \nonumber \P^{(\phi,k,\psi,l)}(T \ge n R  )&  = \E^{(\phi,k,\psi,l)}\bigl[\mathbf{1}_{\{T \ge n R\}}\bigr] = \E^{(\phi,k,\psi,l)}\bigl[\mathbf{1}_{\{T \ge R\}}\E^{(X_{R},\Lambda(R),Y_{R}, \Psi(R))}[\mathbf{1}_{\{T \ge (n-1) R\}}]\bigr]\\
   & \le \E^{(\phi,k,\psi,l)}\bigl[\mathbf{1}_{\{T \ge R\}} \rho^{n-1} \big] \le \rho^n.
\end{align} Obviously, \eqref{induction} also holds for $n=0$. Hence, denoting $\Xi: = \frac{T}{R}$,  we have 
\begin{align}\label{induction2}
 \P^{(\phi,k,\psi,l)}\bigl(\Xi \ge t    \bigr)\leq \P^{(\phi,k,\psi,l)}\bigl(\Xi \ge \lfloor t\rfloor    \bigr) \le \rho^{\lfloor t\rfloor}, \quad \text{ for all } t\geq 0.
\end{align} 
\comment{The proof proceeds by induction on $[t]$: Assume the inequality holds when $t$ is replaced by $[t]- 1$. Then
\beq{anyt}
\begin{split}  \bP^{(\phi,k,\psi,l)}\Bigl(\frac{T}{M+r+N}\geq t\Bigr)
	&\le \mathbb E^{(\phi,k,\psi,l)}\bigl[\mathbf{1}_{\{\frac{T}{M+r+N}\geq[t]\}}\bigr]\\
	&= \mathbb E^{(\phi,k,\psi,l)}\bigl[\mathbf{1}_{\{\frac{T}{M+r+N}\geq[t]\}}\mathbf{1}_{\{\frac{T}{M+r+N}\geq1\}}\bigr]\\
	&= \mathbb E^{(\phi,k,\psi,l)}\bigl[\mathbf{1}_{\{\frac{T}{M+r+N}\geq1\}}\mathbb E^{(X_{1},\Lambda(1),Y_{1},\Lambda^{\prime}(1))}[\mathbf{1}_{\{\frac{T}{M+r+N}\geq[t]-1\}}]\bigr]\\
	&\le \bigl(1-\frac{\delta_2}{2}\bigr)^{[t]-1}\P^{{(\phi,k,\psi,l)}}\Bigl(\frac{T}{M+r+N}\geq1\Bigr)\\
	& \le \bigl(1-\frac{\delta_2}{2}\bigr)^{[t]},
\end{split}
\eeq
where the third equality follows from the Markov property, the fourth line employs the induction hypothesis, and the final inequality utilizes our uniform estimate from \eqref{kl}. Consequently, we obtain
\beq{expection}
\begin{split}
 \mathbb E^{(\phi,k,\psi,l)}\bigl[T\bigr]&=(M+r+N)\mathbb E^{(\phi,k,\psi,l)}\Bigl[\frac{T}{M+r+N}\Bigr]\\
 &\le (M+r+N)\sum_{n=1}^{\infty}n\bP^{(\phi,k,\psi,l)}\Bigl(n-1\leq \frac{T}{M+r+N}<n\Bigr)\\
 &= (M+r+N)\sum_{n=1}^{\infty}\bP^{(\phi,k,\psi,l)}\Bigl(\frac{T}{M+r+N}\geq n-1\Bigr)\\
 &\le (M+r+N)\sum_{n=1}^{\infty}\bigl(1-\frac{\delta_2}{2}\bigr)^{n-1}=\frac{2(M+r+N)}{\delta_2}.
\end{split}
\eeq

Denote $\Xi: = \frac{T}{M+r+N}$ and $\rho: = 1-\frac{\delta_2}{2}$. We have from \eqref{induction} that $\P^{(\phi,k,\psi,l)}(\Xi \ge t) \le \rho^{\lfloor t\rfloor}$. } 

Then for any $n\in \mathbb N$, we have from \eqref{induction2} that  \begin{align}\label{e-Tn1}
 \nonumber \E^{(\phi,k,\psi,l)}[\Xi^n] & = \int_0^\infty n t^{n-1} \P^{(\phi,k,\psi,l)}(\Xi \ge t) \d t \le \int_0^\infty  n t^{n-1}  \rho^{\lfloor t\rfloor} \d t \\ \nonumber &= \int_0^1 n t^{n-1} \d t + \sum_{k=1}^\infty\int_k^{k+1} nt^{n-1}\rho^{k} \d t\\\nonumber & = 1 + \sum_{k=1}^\infty  \rho^{k} [(k+1)^{n}-k^{n}] \\ \nonumber  &  \le 1+ \sum_{k=1}^\infty  \rho^{k} n (k+1)^{n-1}  \\ 
 & \nonumber  = 1+  \frac{n}{\rho} \Bigg[  \sum_{k=1}^\infty  \rho^{k}  k^{n-1} -\rho \Bigg]\\  & \le \frac{n}{\rho}   \mathrm{Li}_{-(n-1)}(\rho), 
\end{align} where the second  last inequality follows from the mean value theorem, and    $\mathrm{Li}_{-(n-1)}(\rho) = \sum_{k=1}^\infty \rho^k k^{n-1}  $  is the polylogarithms function of negative integer order (see   \url{https://mathworld.wolfram.com/Polylogarithm.html} or    \url{https://dlmf.nist.gov/25.12}). 
 We now claim that 
  \begin{align}\label{e:polyln-bound}
    \mathrm{Li}_{-(n-1)}(\rho) \le \frac{\rho (n-1)!}{(1-\rho)^{n}}, \quad \forall \rho \in [0, 1).
  \end{align}  To see this, we note that 
  \begin{align*}
    \frac{(n-1)!}{(1-\rho)^{n}}  & = (n-1)!+ (n-1)! \sum_{k=1}^\infty  \frac{n(n+1)\dots (n+k-1)}{k!} \rho^k \\ & = (n-1)!+ \sum_{k=1}^\infty  \frac{(n+k-1)(n+k-2)\dots (k+1) k!}{k!} \rho^k \\ &  \ge  (n-1)!+ \sum_{k=1}^\infty (k+1)^{n-1} \rho^k \\ & = (n-1)!+ \frac{1}{\rho}\sum_{k=1}^\infty k^{n-1} \rho^k -1\\ &  \ge \frac{1}{\rho}\sum_{k=1}^\infty k^{n-1} \rho^k = \frac{1}{\rho} \mathrm{Li}_{-(n-1)}(\rho). 
  \end{align*} This of course implies   \eqref{e:polyln-bound}. 
  
  Now, plugging \eqref{e:polyln-bound} into \eqref{e-Tn1} yields 
\comment{ \begin{align*}
  \mathrm{Li}_{-(n-1)}(\rho) = \frac{1}{(1-\rho)^{n}} \sum_{k=0}^{n-1}\Big \langle \begin{matrix} n-1\\ k\end{matrix}\Big\rangle \rho^{n-1-k}, 
\end{align*} where $\Big \langle \begin{matrix} n-1\\ k\end{matrix}\Big\rangle$ is the Eulerian numbers,  which counts the number of permutations of $n-1$ elements with $k$ ascents and is defined as
$$\Big \langle \begin{matrix} n-1\\ k\end{matrix}\Big\rangle =\sum_{j=0}^{k+1}(-1)^j\binom{n}{j} (k-j+1)^{n-1} .$$
 Since $\rho \in [0, 1)$,   we have \begin{align*}
   \sum_{k=0}^{n-1} \Big \langle \begin{matrix} n-1\\ k\end{matrix}\Big\rangle \rho^{n-1-k} < \sum_{k=0}^{n-1} \Big \langle \begin{matrix} n-1\\ k\end{matrix}\Big\rangle= (n-1)!. 
  \end{align*} Plugging this into \eqref{e-Tn1} yields }
   \begin{align*}
  \E^{(\phi,k,\psi,l)}[\Xi^n] &     \le  \frac{n!}{ (1-\rho)^{n}}.
\end{align*} Recall that $\Xi = \frac{T}{R}$.  Thus it follows that 
\begin{align}\label{e:T-moments}
  \E^{(\phi,k,\psi,l)}[T^n] & = R^n \E^{(\phi,k,\psi,l)}[\Xi^n] \le R^n   \frac{n!}{  (1-\rho)^{n}} = n! \wdh R^n,  
\end{align}   where $\wdh R = \frac{R}{1-\rho} = \frac{2(M+r+N)}{\delta_2} = 4 e^{H(M+r)} (M +r + N)$. 

Note that for every $\gamma > 2$, we have     $\lambda: = \lambda(\gamma) = \wdh R^{-1} (1- \frac2\gamma) < \wdh R^{-1}$. Then it follows from \eqref{e:T-moments} that \begin{align}\label{lambdaT}
  \E^{(\phi,k,\psi,l)}[{\rm e}^{\lambda T}] & = \sum_{n=0}^\infty \frac{\lambda^n}{n!} \E^{(\phi,k,\psi,l)}[T^n] \le 
   \sum_{n=0}^\infty  \big(\lambda \wdh R\big)^n  = 
    \frac{1}{ 1-\lambda \wdh R} = \frac\gamma 2< \infty.
\end{align}   Note that the inequality \eqref{lambdaT} holds uniformly for all $(\phi, k, \psi, l)\in\mathcal C\times\ss\times \mathcal C\times \ss$. Then 
\begin{align*}
   \mathbb{P}^{(\phi,k,\psi,l)}(T > t)
	& = \mathbb{P}^{(\phi,k,\psi,l)}(e^{\lambda T} > e^{\lambda t})
	\le e^{-\lambda t} \mathbb{E}^{(\phi,k,\psi,l)}[e^{\lambda T}] \\ & \le e^{-\lambda t} \sup_{(\phi,k,\psi,l)\in \C\times \ss\times \C\times\ss} \mathbb{E}^{(\phi,k,\psi,l)}[e^{\lambda T}] \le \frac\gamma2 e^{-\lambda t}.
\end{align*} Again, this holds uniformly for all  $(\phi, k, \psi, l)\in\mathcal C\times\ss\times \mathcal C\times \ss$. Finally, we use the coupling inequality as in the proof of Corollary \ref{success3} to conclude that \begin{align*}
		\|P(t,(\phi,k), \cdot)-\pi (\cdot) \|_{\TV}\le 2
		\displaystyle\int_{E} \d \pi (\psi,l) \P^{(\phi,k,\psi,l)} (T> t)  \le \gamma e^{-\lambda t}.
	\end{align*} As observed before, this estimate holds uniformly for all $(\phi, k)\in \C\times \ss$. Therefore \eqref{eq:uni-erg} is established. This completes the proof. \qed

\begin{Remark}\label{Rem-5.3}
	The parameter $\gamma$ in \eqref{eq:uni-erg} admits a natural interpretation as an allowable initial discrepancy in the total variation bound:   larger values of $\gamma$ permit a greater initial deviation and, in turn, yield faster exponential decay rates, as reflected  by the definition of $\lambda(\gamma)$ in \eqref{eq:lambda_defn}. This trade-off explains why the uniform moment bound on the coupling time $T$ established in \eqref{lambdaT} suffices to guarantee both strong ergodicity and an explicit quantification of the convergence rate. 
\end{Remark}

Combining Theorem~\ref{thm:stroergo} with \cite[Theorem 2.2]{Mao06}, we obtain the following upper bounds for the convergence rates $\beta(\gamma)$ in strong ergodicity of the Markov process $(X_t,\La(t))$.
\begin{Corollary}
	In the setting of Theorem \ref{thm:stroergo}, we have
	\begin{align*}
		\beta(\gamma)\leq \inf\Bigg\{  \left[\frac{2}{\pi(A)}\log\frac{\gamma}{\pi(A)}\right]\left(\sup_{(\phi, k)\in E}\mathbb E^{(\phi, k)}[\tau_A]\right)^{-1}\!: \ A\in \mathcal B(E) \text{ is closed with }\pi(A)> 0\Bigg\},
	\end{align*}
	where $\tau_A$ denotes the hitting time of $A$, i.e., $\tau_A=\inf\{t\geq0: (X_t,\La(t))\in A\}$.
\end{Corollary}

\begin{Remark}\label{Rem-comp-lit}
Theorems~\ref{success1} and~\ref{thm:stroergo} establish that, under Assumptions~\ref{I1}, \ref{assump-Q}, and \ref{Tk}, there exists a successful coupling for the segment process $(X_t,\La(t))$, and that $(X_t,\La(t))$ is strongly ergodic with an explicit exponential convergence rate that is uniformly bounded away from zero. Our results therefore extend the corresponding results for switching diffusions without delay. To the best of our knowledge, this work provides the first exploration of successful coupling and strong ergodicity for switching diffusions with past-dependent switching. 

Even in the delay-free setting (so the transition rates $q_{kl}$ of \eqref{eq:Lambda} depends only on $\phi(0)$), 
our assumptions are significantly weaker than those imposed in the existing literature. Specifically, the models considered in \cite{Xi-13,ShaoX-13,XiShao-13} all require the discrete component $\La$ to have a finite state space. Moreover, the switching component in \cite{Xi-13,ShaoX-13} is a continuous-time Markov chain, while \cite{ShaoX-13} further assumes that the diagonal $\triangle$ is absorbing for the coupled process.

In contrast, our framework allows $\La$ to take values in a countable state space and permits its transition rates to depend on the past trajectory of the continuous component. Furthermore, Assumption~\ref{Tk}(ii) only requires that, for the basic coupling $(\La,\Psi)$ defined by \eqref{eq-Q(x)-coupling}, the rate of entering the diagonal $\triangle$ from any off-diagonal state is uniformly bounded from below. This condition is substantially easier to verify than the classical Doeblin condition, even for continuous-time Markov chains. Finally, we note that Assumptions~\ref{I1}, \ref{assump-Q}, and \ref{Tk}(i) are fairly standard and are consistent with assumptions commonly imposed in the existing literature; see, e.g., the aforementioned references.
\end{Remark}

\section{Application: \texorpdfstring{$N$}{N}-Body Mean Field Model with  Past-State-Dependent Switching}\label{sect-examples}
In this section, we consider an $N$-body mean field model with a past-state-dependent switching process; it is a generalization of the model considered in \cite{XiYin-09}. The model is described by the following SDE: 
\begin{align}\label{e:mean-field}
  \d X_i(t) = \big[\alpha(\La(t))  X_i(t) -X_i^3(t) -\beta(\La(t)) (X_i(t)-\lbar X(t))\big] \d t  +    \sigma_i(X(t), \La(t))  \d W_i(t), \end{align} for $i =1,\dots, N$, 
where $X(t) = (X_1(t),  \dots, X_N(t))^\top\in \R^N$ is the state vector of the system,  $\lbar X(t) = \frac{1}{N}\sum_{i=1}^N X_i(t)$ is the arithmetic mean of the ensemble, $\La(t)$ is a past-state-dependent switching process taking values in $\mathbb S$ and satisfies \eqref{eq:Lambda},
  and $W_i(t), i=1, \dots, N$ are independent standard one-dimensional Brownian motions.  In \eqref{e:mean-field},    $\alpha( k) $ and $\beta(k)  $  are positive constants for each $k\in \ss$, and the coefficient $\sigma_i : \R^N\times \ss \mapsto \R$ is a Borel-measurable function for each $i =1, \dots, N$. 

 As introduced earlier, the model \eqref{e:mean-field} generalizes the framework studied in \cite{XiYin-09}. In their work, the switching component $\Lambda(\cdot)$ has a finite state space, and the transition rates $q_{kl}(\cdot)$ depend solely on the current state of $X$. In contrast, our model permits $\Lambda(\cdot)$ to take values in a countable state space, with transition rates depending on  the past trajectory of    $X$. Furthermore, we relax several restrictive technical assumptions on $q_{kl}(\cdot)$ and $\sigma_i(\cdot)$ imposed in \cite{XiYin-09}.

Notably, \cite{XiYin-09} established strong ergodicity for the process $(X(t), \Lambda(t))$ under additional constraints—namely, that $q_{kl}(x) \equiv q_{kl} > 0$ for all $k\neq l$ and that $\sigma_1(x,k) = \sigma_2(x,k) = \dots = \sigma_N(x,k)$. In Theorem \ref{thm:mf-strongergo} below, we prove strong ergodicity for $(X_t, \Lambda(t))$ without these restrictions. These extensions and relaxations result in a more flexible and comprehensive framework for analyzing mean-field models with path-dependent switching.

For the convenience of presentation, for $(x, k)\in \R^N\times \ss$ and $i =1,\dots, N$, we set 
\begin{align*}
  b_i(x, k) & = \alpha(k) x_i - x_i^3 - \beta(k) (x_i - \lbar x), 
\end{align*}  where $\lbar x = \frac1N\sum_{j=1}^N x_j$,  and 
\begin{align*} b(x, k) & = \big(b_1(x, k), b_2(x, k), \ldots, b_N(x, k)\big)^\top, \\ \sigma(x, k) & = \mathrm{diag}\big(\sigma_1(x, k), \sigma_2(x, k), \ldots, \sigma_N(x, k)\big).
\end{align*} Then   \eqref{e:mean-field} can be rewritten as
\begin{align}\label{e:mean-field2}
  \d X(t) = b(X(t), \La(t)) \d t + \sigma(X(t), \La(t)) \d W(t),
\end{align} where $W(t) = (W_1(t), W_2(t), \ldots, W_N(t))^\top$ is the $N$-dimensional standard Brownian motion. 


\begin{Lemma}\label{lem:mf-existence}
  Assume \eqref{(1.6)} holds and   \begin{align}\label{eq:al-bdd}
    \sup\{\alpha(k): k\in\ss\} & < \infty,  \\ 
     \label{e:sg-growth}|\sigma(x,k)|^2 & \le C(1+|x|^2), \quad \forall (x,k)\in\R^N\times\ss,\\
 \label{e:sg-Lip} |\sigma(x,k)-\sigma(y,k)| & \le C|x-y|, \quad \forall (x,y)\in\R^N\times\R^N, \ k\in\ss,
  \end{align} where $C$ is a positive constant. Then there exists a unique strong solution $(X(t), \La(t))$ to the system \eqref{e:mean-field2} and \eqref{eq:Lambda}.
\end{Lemma}

\para{Proof.} We need to verify that Assumption \ref{I1} holds. Apparently, the coefficients $b(x, k)$ and $\sigma(x, k)$ satisfy the local Lipschitz condition in \eqref{e:b/sig-Lip}. It remains to show that Assumption (ii) is satisfied. To this end, we consider the function $V(x,k) : = |x|^2+1, (x, k)\in \R^N\times \ss$ and note that  \begin{align}\label{e:basic-inequality}
  \sum_{i=1}^N x_i^2 - \lbar x\sum_{i=1}^{N}x_i = \sum_{i=1}^N x_i^2- \frac1 N\Bigg(\sum_{i=1}^{N}x_i\Bigg)^2  \ge 0. 
\end{align} Use this observation in the first inequality below,   it follows  that  \begin{align*}
   \mathcal L_kV (x,k) & =\frac12\tr(\sigma(x,k)\sigma(x,k)^\top 2 I) + \lan 2x, b(x,k)\ran \\ & =  |\sigma(x, k)|^2 +  2 \sum_{i=1}^N  x_i\big(\alpha(k) x_i - x_i^3 - \beta(k) (x_i - \lbar x)\big)   \\ & \le   |\sigma(x, k)|^2 + 2\alpha(k)  |x|^2 - 2\sum_{i=1}^{N} x_i^4 \\ 
   & \le   (C+2\alpha(k))(1+|x|^2) - 2\sum_{i=1}^{N} x_i^4 \\ 
   & = V(x,k) \Bigg[ C+2\alpha(k) -  \frac{2\sum_{i=1}^{N} x_i^4}{ 1+|x|^2}\Bigg] \\& \le -V(x,k) + K, \quad \forall (x,k)\in\R^N\times\ss,
\end{align*} where $K$ is a positive constant, and  the second and third inequalities   above follows from \eqref{e:sg-growth} and \eqref{eq:al-bdd}, respectively.  Thus, Assumption \ref{I1} is satisfied. The existence and uniqueness of the strong solution to the system \eqref{e:mean-field2} and \eqref{eq:Lambda} then follows.  \qed

\begin{Theorem}\label{thm:mf-strongergo}
  Assume that the conditions of Lemma \ref{lem:mf-existence} hold. In addition, suppose that Assumption \ref{assump-Q} and \eqref{qh} are satisfied and that  there exists a positive constant $\lambda_0 \in(0, 1]$ such that \begin{align}\label{e:elliptic}
     \lambda_0 |y|^2 \le \lan  \sigma(x,k)\sigma^\top(x,k) y, y\ran & \le \frac1{ \lambda_0} |y|^2,  
  \end{align}  for all $x,y\in\R^N\times \R^N$   and $k \in \ss$.
  Then the Markov process $(X_t, \La(t))$ is strongly ergodic. 
\end{Theorem} 
\para{Proof.} In view of Theorem \ref{thm:stroergo}, we just need to   verify that Assumption \ref{Tk} (i) holds.  In other words, we need to show that there exists a positive constant $M$ such that \eqref{supTk} holds for all $(x, y, k) \in\R^N\times \R^N\times \ss$.  

Let $\lambda \in (0, \lambda_0)$ be a fixed constant and set $$\sigma_\la(x,k):= \diag(\sigma_{\la,1}(x,k), \dots, \sigma_{\la, N}(x,k)),$$ where $\sigma_{\la, i}(x,k) = \sqrt{\sigma_i(x,k)^2 - \lambda}$ for $(x,k)\in\R^N\times \ss$ and $i =1, \dots, N$.  Thanks to \eqref{e:elliptic}, we have   $\sigma_{\la, i}(x,k) \ge \sqrt { \la_0-\la}$ for all $(x,k)\in \R^N\times \ss $ and $i =1, \dots, N$. Moreover, using \eqref{e:elliptic} and \eqref{e:sg-Lip}, we can derive    \begin{align*}
  |\sigma_{\la,i}(x,k)- \sigma_{\la,i}(y,k)| &  = \frac{|\sigma_i^2(x,k)- \sigma_i^2(y,k)|}{(\sigma_i(x,k)^2 - \lambda )^{\frac12} + (\sigma_i(y,k)^2 - \lambda )^{\frac12}} \\ & \le  \frac{(|\sigma_i(x,k) |+ |\sigma_i(y,k) |)|\sigma_i(x,k) - \sigma_i(y,k)|}{2\sqrt{\la_0-\la}} 
\\ &
  \le  \frac{C|x-y|}{\sqrt{\la_0(\la_0-\la)}}, \quad \forall (x,y)\in\R^N\times\R^N, \ k\in\ss, \ i =1,\dots, N.
\end{align*}  This, of course, implies that  \begin{align}\label{e:sgla-lip} |\sigma_\lambda(x,k)-\sg_\lambda(z,k)|  \le K |x-z|, \qquad  \forall (x,z, k)\in\R^N\times \R^N\times\ss, \end{align} where $K$ is a  positive constant.  

We now  fix an arbitrary $  k \in\ss$ and consider the coupling  process $(X^{(k)}(t), Y^{(k)}(t))$ starting from $(x, y)\in \R^N\times \R^N$ and satisfying the following SDEs:
\begin{align}\label{e:coupling}\begin{cases}
  \d X^{(k)}(t)  = b(X^{(k)}(t), k) \d t + \sigma_\lambda(X^{(k)}(t), k) \d W(t) + \sqrt{\lambda}\, \d B(t), \\ \d Y^{(k)}(t)  = b(Y^{(k)}(t), k) \d t + \sigma_\lambda(Y^{(k)}(t), k) \d W(t) + \sqrt{\lambda} H(X^{(k)}(t), Y^{(k)}(t))
  \d B(t),   
\end{cases}\end{align} where $W(\cdot)$ and $B(\cdot) $ are independent standard $d$-dimensional Brownian motions, and  $$H(x,y) =I-2 \frac{(x-y)(x-y)^\top}{|x-y|^2} \1_{\{x\neq y\}}.$$   Denote the coupling time by  $T^{(k)} = \inf\{t\geq 0: X^{(k)}(t) = Y^{(k)}(t)\}$. As explained in \cite{Wang-23}, the coupling process $(X^{(k)}(t), Y^{(k)}(t))$ given in \eqref{e:coupling} is well-defined.  In addition, we have $X^{(k)}(t) = Y^{(k)}(t)$ for all $t\geq T^{(k)}$. Using  L\'evy's characterization of Brownian motion and the fact that $\lambda I + \sigma^2_\lambda(x,k) = \sigma(x,k)\sigma^\top(x,k)$, we can verify that the marginal distribution of \eqref{e:coupling} agrees with that of the solution of $\d X^{(k)}(t)  = b(X^{(k)}(t), k) \d t + \sigma(X^{(k)}(t), k) \d W(t)$.  
Thus, the process $(X^{(k)}(t), Y^{(k)}(t))$ is indeed a coupling of the   process $X^{(k)}(t) $.  

The operator  $\wdh\Omega_k$ associated with the coupling process $(X^{(k)}(t), Y^{(k)}(t))$ is given by the following. First, we denote \begin{align*}
  a(x, y, k): = \begin{pmatrix} a(x,k) & c(x,y,k)\\ c(x,y,k)^\top & a(y,k) \end{pmatrix}, \text{ and } b(x, y, k): = \begin{pmatrix} b(x,k) \\ b(y,k) \end{pmatrix},
\end{align*} where $c(x,y, k): =  \lambda H(x,y) + \sigma_\lambda(x,k) \sg^\top_\lambda(y, k) $.
One can verify that $a(x, y, k)$ is symmetric and positive definite.  For any $f \in C^2(\R^N\times\R^N)$, we define the operator $\wdh\Omega_k$ by
\begin{align}\label{Omega}
  \wdh\Omega_k f(x,y) & = \frac12 \tr\big(a(x,y,k)D^2\phi(x,y)\big) + \lan b(x,y,k), D\phi(x,y)\ran.
\end{align}  
In particular, for any  $\phi \in C^2([0,\infty))$ and all $x, z \in \mathbb{R}^N$ with $x \neq z$, we have from straightforward calculations that
 \begin{align}\label{Omega_d}
 \wdh\Omega_{k}  \phi(|x-z|) &= \frac12  \phi''(|x-z|)\lbar{A}(x,z,k) 
  + \frac{ \phi' (|x-z|)}{2|x-z|}[\tr A(x,z,k)
  - \lbar{A}(x,z,k) +2B(x,z,k)], 
 \end{align} where \begin{align*}
			A(x,z,k) &:= a(x,k) + a(z,k) - 2c(x,z,k),\\
			\lbar{A}(x,z,k) &:= \frac{1}{|x-z|^2}\langle x-z , A(x,z,k)(x-z)\rangle,\\
			B(x,z,k) &:= \langle x-z, b(x,k) - b(z,k)\rangle. 
		\end{align*} Direct calculations show that \begin{align}
      \label{e:A-ineq}
    \tr A(x,z,k) = |\sigma_{\lambda}(x,k) - \sigma_{\lambda}(z,k)|^2 + 4\lambda,  
     \end{align}  and 
     {\begin{align*}
      \lbar{A} (x,z,k) = 4\lambda + \frac{1}{|x-z|^2} |( \sigma_\lambda(x, k) - \sigma_\lambda(z, k)) (x-z)|^2
     \end{align*} for all $x, z \in \R^N$ with $x\neq z$ and $k\in \ss$.
 Since $\sigma(x,k)\sigma^\top(x,k) = \lambda I + \sigma^2_\lambda(x,k)$, it follows from \eqref{e:elliptic}    that $|\sigma_\lambda(x,k) y|^2 \le (\frac{1}{\lambda_0} -\lambda) |y|^2$ for all $x,y\in\R^N$ and $k\in\ss$. This, in turn, implies that $|\sigma_\lambda(x,k) |^2  \le \theta$ for all $(x, k)\in \R^N\times \ss$, where $\theta$ is a positive constant. Then, it follows that for all $x\neq z \in \R^N$ and $k\in \ss$, we have  \begin{align}\label{e2:A-ineq}
  4\lambda \le  \lbar{A} (x,z,k) \le 4(\lambda + \theta), \quad \forall (x,z,k)\in\R^N\times\R^N\times\ss \text{ with }x\neq z.
\end{align}}

Next we compute \begin{align*}
  B(x,z,k) & = \langle x-z,b(x,k)-b(z,k)\rangle \\ & = \sum_{i=1}^N (x_i-z_i)\big[ \alpha(k) (x_i-z_i) - (x_i^3-z_i^3) - \beta(k) \big( (x_i-\lbar x) - (z_i-\lbar z)\big )\big]  
   \\   & = \sum_{i=1}^N  \alpha(k) (x_i-z_i)^2 - \sum_{i=1}^{N}  (x_i-z_i)^2(x_i^2 + x_iz_i +z_i^2) \\   & \quad - \beta(k) \sum_{i=1}^{N} (x_i-z_i)^2 + \beta(k) \sum_{i=1}^{N} (x_i-z_i)(\lbar x-\lbar z) \\ 
  & \le \alpha(k) |x-z|^2 - \frac14\sum_{i=1}^{N}(x_i-z_i)^4, \end{align*} where we used \eqref{e:basic-inequality} and the elementary inequality that $x_i^2+x_iz_i + z_i^2 \ge \frac14(x_i-z_i)^2$ to derive the last inequality. Another application of \eqref{e:basic-inequality} gives us \begin{align}\label{e:B-ineq}
    B(x,z,k) & \le \alpha(k) |x-z|^2 - \frac14\sum_{i=1}^{N}(x_i-z_i)^4 \le \alpha(k) |x-z|^2 - \frac1{4N} |x-z|^4.
  \end{align} A combination of \eqref{e:sgla-lip} and \eqref{e:B-ineq} then leads to  
\begin{equation}
\label{eq:coeff-cts}
\begin{aligned}
 |\sigma_{\lambda}(x,k)  - \sigma_{\lambda}(z,k)|^2  + 2 B(x,z,k) 
  & \leq \big(K{^2}+ 2 \alpha(k)\big) |x-z|^2 - \frac1{4N} |x-z|^4 \\ & \le \kappa |x-z|^2 - \frac1{4N} |x-z|^4, 
\end{aligned}
\end{equation} for all $ (x,z, k)\in\R^N\times \R^N\times\ss$, 
  where $\kappa = K^2 + 2\sup_{k\in \ss}\alpha(k)$, which is finite thanks to \eqref{eq:al-bdd}.
 {From \eqref{e:A-ineq}, \eqref{e2:A-ineq}, and \eqref{eq:coeff-cts},  we obtain
 \begin{align}\label{eq:g} 
 		\nonumber \frac{\tr A(x,z,k)-\lbar{A}(x,z,k)+2 B(x,z,k)}{\lbar{A}(x,z,k)} &\leq \frac{\kappa |x-z|^2 - \frac1{4N} |x-z|^4 }{\lbar{A}(x,z,k) } \\ \nonumber  & \le \frac{\kappa}{4\lambda} |x-z|^2 - \frac{1}{16N(\lambda +\theta)} |x-z|^4\\
 		&=|x-z|g(|x-z|),
 	\end{align}
where $$g(r) = \frac{\kappa}{4\lambda} r- \frac{1}{16N(\lambda +\theta)} r^3, \quad  \text{  for }r\geq 0.$$ }

  Motivated by
  \cite{PriolaW-06},  
   we now consider the function $G$  defined by
 \begin{align*}
   G(\rho):= \int_{0}^{\rho} f(s) 
   \d s, \ \ \rho \ge 0, 
 \end{align*} where $$f(s): =\exp\bigg\{-\int_{0}^{s} g(w) \d w\bigg\} \int_{s}^{\infty}\exp \bigg\{\int_{0}^{v} g(u) \d u\bigg\}\d v, \ \ s \ge 0.  $$  Note that   $G(\rho)$ is well-defined for each $\rho \ge 0$ because   $$\int_{s}^{\infty}\exp \bigg\{\int_{0}^{v} g(u) \d u\bigg\}\d v   = \int_{s}^{\infty}\exp \bigg\{\frac{\kappa}{8\lambda} v^2- \frac{1}{64N(\lambda+\theta)} v^4\bigg\}\d v< \infty,$$ for each $s > 0$. Obviously, we have  $G(0)= \lim_{\rho\rightarrow 0} G(\rho) = 0$. In addition, direct calculations show that the function $G$ is twice continuously differentiable with 
\begin{equation}
\label{e:G'>0}
\begin{aligned}
&   G'(\rho) = f(\rho) \geq 0,
 \text{ and}\ \
G''(\rho) = -1 - g(\rho)  G'(\rho).
 \end{aligned}
\end{equation} Hence    $G(\rho) \ge 0$ for all $\rho \ge 0$.

 We  next verify that   $G(\infty): = \lim_{\rho \to \infty} G(\rho)   < \infty$. 
 To this end, we consider the auxiliary   function 
 \begin{align*}
  h(s):=\frac{16N(\lambda+\theta)}{ s^3}, \quad s >0. 
 \end{align*}
 Obviously, we have 
 \begin{equation}\label{eq:h}
 	   \int_{1}^{\infty}h(s)\d s<\infty.
 \end{equation}
 Thanks to  l'H\^opital's rule and   straightforward calculations, we can show  that $\lim_{s \to \infty} \frac{f(s)}{h(s)} = 1$.
 Combining this result with \eqref{eq:h} yields $\int_1^\infty f(s) \d s < \infty$.
 Note also that   $\int_0^1 f(s) \d s < \infty$. Therefore  it follows that 
 \[
 G(\infty) = \int_0^\infty f(s) \d s < \infty.
 \]
  Consequently, $G$ is nonnegative and  uniformly bounded.

We now use \eqref{Omega_d}, \eqref{e:A-ineq}, \eqref{eq:g}, and \eqref{e:G'>0} to  compute 
\begin{align}\label{e:Omega-G-estimate}
  \nonumber \wdh\Omega_{k}  G(|x-z|) & =  \frac{  G''(|x-z|)}{2}\lbar{A}(x,z,k) + \frac{G'(|x-z|)}{2|x-z|}\big[\tr A(x,z,k) - \lbar{A}(x,z,k) + 2B(x,z,k)\big] \\ 
   \nonumber  & = \Bigg[-\frac12 - \frac{1}{2}g(|x-z|)  G'(|x-z|)  \\ 
    \nonumber  & \qquad + \frac{G'(|x-z|)}{2|x-z|}\frac{  \tr A(x,z,k) - \bar{A}(x,z,k) + 2B(x,z,k)}{\lbar{A}(x,z,k) }\Bigg]\lbar{A}(x,z,k)   \\
  \nonumber   & \le  \Bigg[-\frac12 - \frac{1}{2}g(|x-z|)  G'(|x-z|)  + \frac{G'(|x-z|)}{2|x-z|} |x-z|g(|x-z|)\Bigg]\lbar{A}(x,z,k) \\
  & = -\frac{ \lbar{A}(x,z,k) }{2} \le -2\lambda,  
\end{align} for all $(x,z,k)\in\R^N\times\R^N\times\ss $ with $ |x-z|> 0$.  

Finally, for the coupling process $(X^{(k)}(t), Y^{(k)}(t))$ given in \eqref{e:coupling} with initial condition $(x, y)\in \R^N\times \R^N$ and $x\neq y$, let $\beta_n : = \inf\{t\ge 0: |X^{(k)}(t)|\vee |Y^{(k)}(t)| \ge n \}$, $n\in \mathbb N$. For any $t \ge 0$, it follows from   It\^o's formula, \eqref{e:Omega-G-estimate},  and the optional stopping theorem that \begin{align*}
  \mathbb E& ^{(x,y,k)}\big[G(|X^{(k)}(T^{(k)}\wedge \beta_n\wedge t)-Y^{(k)}(T^{(k)}\wedge \beta_n\wedge t)|)\big]
  \\ & = G(|x-y|) +  \mathbb E^{(x,y,k)}\bigg[\int_0^{T^{(k)}\wedge \beta_n\wedge t} \wdh\Omega_{k} G(|X^{(k)}(s)-Y^{(k)}(s)|) \d s\bigg]    \\ & \le G(|x-y|)   -2\lambda \mathbb E^{(x,y,k)}[T^{(k)}\wedge \beta_n\wedge t].
\end{align*} Since $G$ is nonnegative and uniformly bounded, we can take the limit as $n\to \infty$ and then as $t\to \infty$ to obtain $$ \mathbb E^{(x,y,k)}[T^{(k)}] \le \frac{G(|x-y|)}{2\lambda} \le K,$$ where $K$ is a positive constant independent of $x, y$, and $k$.   This of course implies that Assumption \ref{Tk} (i) and hence the process $(X_t, \La(t))$ is strongly ergodic thanks to Theorem \ref{thm:stroergo}.  The proof is complete. \qed



\begin{thebibliography}{}\parskip=-2pt

\bibitem[Applebaum, 2009]{APPLEBAUM}
Applebaum, D. (2009).
\newblock {\em L\'evy processes and stochastic calculus}, volume 116 of {\em Cambridge Studies in Advanced Mathematics}.
\newblock Cambridge University Press, Cambridge, second edition.

\bibitem[Bao and Shao, 2016]{BaoS-16}
Bao, J. and Shao, J. (2016).
\newblock Permanence and extinction of regime-switching predator-prey models.
\newblock {\em SIAM J. Math. Anal.}, 48(1):725--739.

\bibitem[Cai et~al., 2021]{CCM-21}
Cai, S., Cai, Y., and Mao, X. (2021).
\newblock A stochastic differential equation {SIS} epidemic model with regime switching.
\newblock {\em Discrete Contin. Dyn. Syst. Ser. B}, 26(9):4887--4905.

\bibitem[Cao et~al., 2024]{CWW-24}
Cao, W., Wu, F., and Wu, M. (2024).
\newblock Weak convergence and stability of functional diffusion systems with singularly perturbed regime switching.
\newblock {\em Nonlinear Anal. Hybrid Syst.}, 53:Paper No. 101487, 17.

\bibitem[Chen, 2004]{Chen04}
Chen, M.-F. (2004).
\newblock {\em From {M}arkov chains to non-equilibrium particle systems}.
\newblock World Scientific Publishing Co. Inc., River Edge, NJ, second edition.

\bibitem[Chen and Li, 1989]{ChenLi-89}
Chen, M.~F. and Li, S.~F. (1989).
\newblock Coupling methods for multidimensional diffusion processes.
\newblock {\em Ann. Probab.}, 17(1):151--177.

\bibitem[Da~Prato and Zabczyk, 1996]{DaPratoZ-96}
Da~Prato, G. and Zabczyk, J. (1996).
\newblock {\em Ergodicity for infinite-dimensional systems}, volume 229 of {\em London Mathematical Society Lecture Note Series}.
\newblock Cambridge University Press, Cambridge.

\bibitem[Da~Prato and Zabczyk, 2014]{DaPratoZ-14}
Da~Prato, G. and Zabczyk, J. (2014).
\newblock {\em Stochastic equations in infinite dimensions}, volume 152.
\newblock Cambridge university press.

\bibitem[Elliott and Siu, 2010]{Elliott-Siu-10}
Elliott, R.~J. and Siu, T.~K. (2010).
\newblock On risk minimizing portfolios under a {M}arkovian regime-switching {B}lack-{S}choles economy.
\newblock {\em Ann. Oper. Res.}, 176:271--291.

\bibitem[Greenhalgh et~al., 2016]{GLM-16}
Greenhalgh, D., Liang, Y., and Mao, X. (2016).
\newblock Modelling the effect of telegraph noise in the {SIRS} epidemic model using {M}arkovian switching.
\newblock {\em Phys. A}, 462:684--704.

\bibitem[Huang et~al., 2016]{Huang-16}
Huang, J., Zhang, H., and Zhang, J. (2016).
\newblock A unified approach to diffusion analysis of queues with general patience-time distributions.
\newblock {\em Math. Oper. Res.}, 41(3):1135--1160.

\bibitem[Li et~al., 2017]{LLC-17}
Li, D., Liu, S., and Cui, J. (2017).
\newblock Threshold dynamics and ergodicity of an {SIRS} epidemic model with {M}arkovian switching.
\newblock {\em J. Differential Equations}, 263(12):8873--8915.

\bibitem[Li et~al., 2022]{LLLM-22}
Li, X., Liu, W., Luo, Q., and Mao, X. (2022).
\newblock Stabilisation in distribution of hybrid stochastic differential equations by feedback control based on discrete-time state observations.
\newblock {\em Automatica J. IFAC}, 140:Paper No. 110210, 8.

\bibitem[Li et~al., 2020]{LMMY-20}
Li, X., Mao, X., Mukama, D.~S., and Yuan, C. (2020).
\newblock Delay feedback control for switching diffusion systems based on discrete-time observations.
\newblock {\em SIAM J. Control Optim.}, 58(5):2900--2926.

\bibitem[Lindvall and Rogers, 1986]{LindR-86}
Lindvall, T. and Rogers, L. C.~G. (1986).
\newblock Coupling of multidimensional diffusions by reflection.
\newblock {\em Ann. Probab.}, 14(3):860--872.

\bibitem[Liu et~al., 2012]{LiuMW:12}
Liu, X., Margaritis, D., and Wang, P. (2012).
\newblock Stock market volatility and equity returns: Evidence from a two-state markov-switching model with regressors.
\newblock {\em Journal of Empirical Finance}, 19(4):483--496.

\bibitem[Luo and Mao, 2009]{Luo-Mao09}
Luo, Q. and Mao, X. (2009).
\newblock Stochastic population dynamics under regime switching. {II}.
\newblock {\em J. Math. Anal. Appl.}, 355(2):577--593.

\bibitem[Mao and Yuan, 2006]{MaoY}
Mao, X. and Yuan, C. (2006).
\newblock {\em Stochastic differential equations with {M}arkovian switching}.
\newblock Imperial College Press, London.

\bibitem[Mao, 2002]{Mao02}
Mao, Y.-H. (2002).
\newblock Strong ergodicity for {M}arkov processes by coupling methods.
\newblock {\em J. Appl. Probab.}, 39(4):839--852.

\bibitem[Mao, 2006]{Mao06}
Mao, Y.-H. (2006).
\newblock Convergence rates in strong ergodicity for {M}arkov processes.
\newblock {\em Stochastic Process. Appl.}, 116(12):1964--1976.

\bibitem[Meyn and Tweedie, 2009]{MeynT-09}
Meyn, S. and Tweedie, R.~L. (2009).
\newblock {\em Markov chains and stochastic stability}.
\newblock Cambridge University Press, Cambridge, second edition.
\newblock With a prologue by Peter W. Glynn.

\bibitem[Nguyen et~al., 2021]{NNN-21}
Nguyen, D.~H., Nguyen, D., and Nguyen, S.~L. (2021).
\newblock Stability in distribution of path-dependent hybrid diffusion.
\newblock {\em SIAM J. Control Optim.}, 59(1):434--463.

\bibitem[Nguyen et~al., 2024]{NNTY-23}
Nguyen, D.~H., Nguyen, N.~N., Ta, T.~T., and Yin, G. (2024).
\newblock Stability of stochastic functional differential equations with past-dependent random switching involving countably infinite states.
\newblock {\em IEEE Trans. Automat. Control}, 69(3):1612--1626.

\bibitem[Nguyen and Yin, 2016]{NguyenY-16}
Nguyen, D.~H. and Yin, G. (2016).
\newblock Modeling and analysis of switching diffusion systems: past-dependent switching with a countable state space.
\newblock {\em SIAM J. Control Optim.}, 54(5):2450--2477.

\bibitem[Nguyen and Yin, 2020]{NguyenY-20}
Nguyen, D.~H. and Yin, G. (2020).
\newblock Stability of stochastic functional differential equations with regime-switching: analysis using {D}upire's functional {I}t\^{o} formula.
\newblock {\em Potential Anal.}, 53(1):247--265.

\bibitem[Priola and Wang, 2006]{PriolaW-06}
Priola, E. and Wang, F.-Y. (2006).
\newblock Gradient estimates for diffusion semigroups with singular coefficients.
\newblock {\em J. Funct. Anal.}, 236(1):244--264.

\bibitem[Savku, 2026]{Savku-24}
Savku, E. (2026).
\newblock An approach for regime-switching stochastic control problems with memory and terminal conditions.
\newblock {\em Optimization}, 75(2):395--412.

\bibitem[Schaller and Norden, 1997]{SchavN:97}
Schaller, H. and Norden, S.~V. (1997).
\newblock Regime switching in stock market returns.
\newblock {\em Applied Financial Economics}, 7(2):177--191.

\bibitem[Settati and Lahrouz, 2014]{SetLa-14}
Settati, A. and Lahrouz, A. (2014).
\newblock Stationary distribution of stochastic population systems under regime switching.
\newblock {\em Appl. Math. Comput.}, 244:235--243.

\bibitem[Shao, 2015]{Shao-15}
Shao, J. (2015).
\newblock Strong solutions and strong {F}eller properties for regime-switching diffusion processes in an infinite state space.
\newblock {\em SIAM J. Control Optim.}, 53(4):2462--2479.

\bibitem[Shao, 2017]{Shao-17}
Shao, J. (2017).
\newblock Stabilization of regime-switching processes by feedback control based on discrete time observations.
\newblock {\em SIAM J. Control Optim.}, 55(2):724--740.

\bibitem[Shao, 2018]{Shao-18}
Shao, J. (2018).
\newblock Invariant measures and {E}uler-{M}aruyama's approximations of state-dependent regime-switching diffusions.
\newblock {\em SIAM J. Control Optim.}, 56(5):3215--3238.

\bibitem[Shao and Xi, 2013]{ShaoX-13}
Shao, J. and Xi, F. (2013).
\newblock Strong ergodicity of the regime-switching diffusion processes.
\newblock {\em Stochastic Process. Appl.}, 123(11):3903--3918.

\bibitem[Shao and Xi, 2019]{ShaoX-19}
Shao, J. and Xi, F. (2019).
\newblock Stabilization of regime-switching processes by feedback control based on discrete time observations {II}: {S}tate-dependent case.
\newblock {\em SIAM J. Control Optim.}, 57(2):1413--1439.

\bibitem[Skorokhod, 1989]{Skorohod-89}
Skorokhod, A.~V. (1989).
\newblock {\em Asymptotic methods in the theory of stochastic differential equations}, volume~78 of {\em Translations of Mathematical Monographs}.
\newblock American Mathematical Society, Providence, RI.
\newblock Translated from the Russian by H. H. McFaden.

\bibitem[Song et~al., 2011]{Song-S-Z}
Song, Q.~S., Stockbridge, R.~H., and Zhu, C. (2011).
\newblock On optimal harvesting problems in random environments.
\newblock {\em SIAM J. Control Optim.}, 49(2):859--889.

\bibitem[Tran et~al., 2022]{TNY-22}
Tran, K.~Q., Nguyen, D.~H., and Yin, G. (2022).
\newblock Stability in distribution and stabilization of switching jump diffusions.
\newblock {\em ESAIM Control Optim. Calc. Var.}, 28:Paper No. 72, 26.

\bibitem[Wang, 2010]{Wang-10}
Wang, F.-Y. (2010).
\newblock Harnack inequalities on manifolds with boundary and applications.
\newblock {\em J. Math. Pures Appl. (9)}, 94(3):304--321.

\bibitem[Wang, 2023]{Wang-23}
Wang, F.-Y. (2023).
\newblock Exponential ergodicity for non-dissipative {M}c{K}ean-{V}lasov {SDE}s.
\newblock {\em Bernoulli}, 29(2):1035--1062.

\bibitem[Wen et~al., 2023]{WLXZ-23}
Wen, J., Li, X., Xiong, J., and Zhang, X. (2023).
\newblock Stochastic linear-quadratic optimal control problems with random coefficients and {M}arkovian regime switching system.
\newblock {\em SIAM J. Control Optim.}, 61(2):949--979.

\bibitem[Whitt, 2002]{Whitt-02}
Whitt, W. (2002).
\newblock {\em Stochastic-process limits: : An introduction to stochastic-process limits and their application to queues}.
\newblock Springer Series in Operations Research. Springer-Verlag, New York.

\bibitem[Xi, 2013]{Xi-13}
Xi, F. (2013).
\newblock Coupling for {M}arkovian switching jump-diffusions.
\newblock {\em Appl. Math. J. Chinese Univ. Ser. B}, 28(2):204--216.

\bibitem[Xi and Shao, 2013]{XiShao-13}
Xi, F. and Shao, J. (2013).
\newblock Successful couplings for diffusion processes with state-dependent switching.
\newblock {\em Sci. China Math.}, 56(10):2135--2144.

\bibitem[Xi and Yin, 2009]{XiYin-09}
Xi, F. and Yin, G. (2009).
\newblock Asymptotic properties of a mean-field model with a continuous-state-dependent switching process.
\newblock {\em J. Appl. Probab.}, 46(1):221--243.

\bibitem[Xi and Yin, 2010]{XiYin-10}
Xi, F. and Yin, G. (2010).
\newblock Asymptotic properties of nonlinear autoregressive {M}arkov processes with state-dependent switching.
\newblock {\em J. Multivariate Anal.}, 101(6):1378--1389.

\bibitem[Xi and Yin, 2011]{XiYin-11}
Xi, F. and Yin, G. (2011).
\newblock Jump-diffusions with state-dependent switching: existence and uniqueness, {F}eller property, linearization, and uniform ergodicity.
\newblock {\em Sci. China Math.}, 54(12):2651--2667.

\bibitem[Xi and Yin, 2015]{XiYin-15}
Xi, F. and Yin, G. (2015).
\newblock Stochastic {L}i{\'e}nard equations with state-dependent switching.
\newblock {\em Acta Math. Appl. Sin. Engl. Ser.}, 31(4):893--908.

\bibitem[Xi and Zhao, 2006]{XiZ-06}
Xi, F. and Zhao, L. (2006).
\newblock On the stability of diffusion processes with state-dependent switching.
\newblock {\em Sci. China Ser. A}, 49(9):1258--1274.

\bibitem[Xi and Zhu, 2017]{XiZ-17}
Xi, F. and Zhu, C. (2017).
\newblock On {F}eller and strong {F}eller properties and exponential ergodicity of regime-switching jump diffusion processes with countable regimes.
\newblock {\em SIAM J. Control Optim.}, 55(3):1789--1818.

\bibitem[Xi et~al., 2021]{XiZW-21}
Xi, F., Zhu, C., and Wu, F. (2021).
\newblock On strong {F}eller property, exponential ergodicity and large deviations principle for stochastic damping {H}amiltonian systems with state-dependent switching.
\newblock {\em J. Differential Equations}, 286:856--891.

\bibitem[Yin and Zhu, 2010]{YZ-10}
Yin, G.~G. and Zhu, C. (2010).
\newblock {\em Hybrid Switching Diffusions: Properties and Applications}, volume~63 of {\em Stochastic Modelling and Applied Probability}.
\newblock Springer, New York.

\bibitem[Zhang, 2001]{Zhang-01}
Zhang, Q. (2001).
\newblock Stock trading: an optimal selling rule.
\newblock {\em SIAM J. Control Optim.}, 40(1):64--87.

\bibitem[Zhou and Yin, 2003]{Zhou-Yin}
Zhou, X.~Y. and Yin, G. (2003).
\newblock Markowitz's mean-variance portfolio selection with regime switching: a continuous-time model.
\newblock {\em SIAM J. Control Optim.}, 42(4):1466--1482.

\bibitem[Zhu and Yin, 2009]{ZY-09c}
Zhu, C. and Yin, G. (2009).
\newblock On hybrid competitive {L}otka-{V}olterra ecosystems.
\newblock {\em Nonlinear Anal.}, 71(12):e1370--e1379.

\bibitem[Zhu et~al., 2015]{ZhuYB-15}
Zhu, C., Yin, G., and Baran, N.~A. (2015).
\newblock Feynman-{K}ac formulas for regime-switching jump diffusions and their applications.
\newblock {\em Stochastics}, 87(6):1000--1032.

\end{thebibliography}
 
\def\cprime{$'$} \def\polhk#1{\setbox0=\hbox{#1}{\ooalign{\hidewidth \lower1.5ex\hbox{`}\hidewidth\crcr\unhbox0}}}

\end{document}